\documentclass[10pt,english]{article}
\usepackage{geometry}
 \usepackage{xcolor}
\definecolor{ured}{RGB}{0, 0, 0}
\definecolor{ublue}{RGB}{0, 90, 255}
\definecolor{ugreen}{RGB}{3, 175, 122}
\usepackage{array}
\usepackage{bbding}
\usepackage{multirow}
\usepackage[fleqn]{amsmath}
\usepackage{amssymb}
\usepackage{amsthm}
\usepackage{graphicx}
\usepackage{natbib}
\usepackage{lscape}
\usepackage[hyphens]{url}
\usepackage[hidelinks]{hyperref}
\hypersetup{
    colorlinks=true,
    citecolor=blue,
    linkcolor=blue,
    urlcolor=black,
}
\usepackage{comment}
\usepackage[title]{appendix}
\usepackage{longtable}
\usepackage{mathtools}
\usepackage{chngcntr}
\usepackage{authblk}
\usepackage{babel}
\usepackage{subcaption}
\usepackage{fancyhdr}
\usepackage{cleveref}
\usepackage{abstract}

\makeatletter
\allowdisplaybreaks
\date{}

\fancypagestyle{plain}{\fancyhf{} 
	\chead{\footnotesize \textcolor{gray}{\textit{\shortauthors -- \runningtitle}}}
	\rhead{\footnotesize \textcolor{gray}{\thepage}}
}
\makeatletter
\def\blfootnote{\xdef\@thefnmark{}\@footnotetext}
\makeatother

\makeatletter
\def\titlepageext{
	\begin{center}	
	{\parindent0pt
		\rule{0.9\textwidth}{1pt}
		\begin{minipage}[t]{0.30\textwidth}
			\small {\it Keywords:}\\
			\keyword
		\end{minipage}\hspace{3mm}
		\begin{minipage}[t]{0.6\textwidth}
			\small \abstract
		\end{minipage}

		\rule{0.9\textwidth}{2pt}
	}\end{center}

	\blfootnote{* Corresponding author. 
  E-mail address: 
  \href{mailto:\corresemail}{sakai.t.dcad@m.isct.ac.jp} (T. Sakai),
  \\
  \href{mailto:\corresemail}{takashi.akamatsu.c4@tohoku.ac.jp} (T. Akamatsu),
  \href{mailto:\corresemail}{satsukawa@tohoku.ac.jp} (K. Satsukawa).
  }
}
\makeatother

\usepackage[mathlines, pagewise]{lineno}
\usepackage{etoolbox} 

\newcommand*\linenomathpatchAMS[1]{\expandafter\pretocmd\csname #1\endcsname {\linenomathAMS}{}{}\expandafter\pretocmd\csname #1*\endcsname{\linenomathAMS}{}{}\expandafter\apptocmd\csname end#1\endcsname {\endlinenomath}{}{}\expandafter\apptocmd\csname end#1*\endcsname{\endlinenomath}{}{}}

\expandafter\ifx\linenomath\linenomathWithnumbers
\let\linenomathAMS\linenomathWithnumbers
\patchcmd\linenomathAMS{\advance\postdisplaypenalty\linenopenalty}{}{}{}
\else
\let\linenomathAMS\linenomathNonumbers
\fi

\linenomathpatchAMS{gather}
\linenomathpatchAMS{multline}
\linenomathpatchAMS{align}
\linenomathpatchAMS{alignat}
\linenomathpatchAMS{flalign}

\nolinenumbers

\usepackage[ISO]{diffcoeff}
\usepackage[utf8]{inputenc}
\usepackage[T1]{fontenc}
\usepackage{newpxtext, newpxmath}

\usepackage{tikz}
\usepackage{pgfplots}
\usetikzlibrary{shapes.geometric}
\usetikzlibrary{arrows.meta}
\usetikzlibrary{overlay-beamer-styles}
\usetikzlibrary{calc,3d}
\usetikzlibrary{decorations.pathreplacing}
\usetikzlibrary{math}
 \usetikzlibrary{patterns}
 
\usepackage{wasysym}
\makeatletter
\newcommand{\UE}[1]{#1^{\mathrm{UE}}} 
\newcommand{\CD}[1]{#1^{\mathrm{C}}}
\newcommand{\FD}[1]{#1^{\mathrm{F}}}
\newcommand{\SO}[1]{#1^{\mathrm{SO}}}
\makeatother

\theoremstyle{definition}
\newtheorem{thm}{Theorem}[section]
\newtheorem{cor}{Corollary}[section]
\newtheorem{dfn}{Definition}[section]
\newtheorem{lem}{Lemma}[section]
\newtheorem{asm}{Assumption}[section]
\newtheorem{pro}{Proposition}[section]

\newtheorem*{prf}{Proof}

\crefname{section}{Section}{Sections}
\crefname{subsection}{Section}{Sections}

\crefname{equation}{Eq.}{Eqs.}
\crefname{figure}{Figure}{Figures}
\crefname{subfigure}{Figure}{Figures}
\crefname{table}{Table}{Tables}

\crefname{thm}{Theorem}{Theorems}
\crefname{cor}{Corollary}{Corollary}
\crefname{dfn}{Definition}{Definition}
\crefname{lem}{Lemma}{Lemmas}
\crefname{asm}{Assumption}{Assumption}
\crefname{pro}{Proposition}{Proposition}
\crefname{cnj}{Conjecture}{Conjecture}

\usepackage{algorithm}
\usepackage{algorithmic}

\usepackage{bm}
\newcommand{\Vt}[1]{\bm{#1}}

\newcommand{\Vtb}{\bm{b}}
\newcommand{\Vtc}{\bm{c}}

\newcommand{\Vtp}{\bm{p}}
\newcommand{\Vtq}{\bm{q}}

\newcommand{\Vts}{\bm{s}}

\newcommand{\Vtw}{\bm{w}}
\newcommand{\Vtx}{\bm{x}}
\newcommand{\Vty}{\bm{y}}

\newcommand{\VtA}{\bm{A}}
\newcommand{\VtB}{\bm{B}}

\newcommand{\VtM}{\bm{M}}

\newcommand{\VtQ}{\bm{Q}}

\newcommand{\Vttheta}{\bm{\theta}}

\newcommand{\Vtlambda}{\bm{\lambda}}
\newcommand{\Vtmu}{\bm{\mu}}

\newcommand{\Vtpi}{\bm{\pi}}

\newcommand{\Vtrho}{\bm{\rho}}

\newcommand{\Vtsigma}{\bm{\sigma}}

\newcommand{\ClL}{\mathcal{L}}

\newcommand{\ClN}{\mathcal{N}}

\newcommand{\ClP}{\mathcal{P}}

\newcommand{\ClT}{\mathcal{T}}

 \usepackage{framed}
\usepackage{tabularx}

\title{Queue Replacement Approach \\ to Dynamic User Equilibrium Assignment \\ with Route and Departure Time Choice}

\def\shortauthors{Sakai et al.}
\def\runningtitle{Queue Replacement Approach to DUE Assignment}

\author[a, $\ast$]{Takara Sakai}
\author[b, $\ast$]{Takashi Akamatsu}
\author[b, $\ast$]{Koki Satsukawa}
\affil[a]{Department of Civil and Environmental Engineering, Institute of Science Tokyo, 2-12-1 W6-9, Ookayama, Meguro, Tokyo 152-8550, Japan}
\affil[b]{Graduate School of Information Sciences, Tohoku University, 6-6 Aramaki Aoba, Aoba-ku, Sendai, Miyagi 980-8579, Japan}

\def\corresemail{sakai.t.dcad@m.isct.ac.jp~(T. Sakai), takashi.akamatsu.c4@tohoku.ac.jp~(T. Akamatsu), satsukawa@tohoku.ac.jp~(K. Satsukawa).}

\def\abstract{
This study develops a hybrid analytical and numerical approach for dynamic user equilibrium (DUE) assignment with simultaneous route and departure time choice (RDTC) for homogeneous users.
The core concept of the proposed approach is the generalized queue replacement principle (GQRP), which establishes an equivalence between the equilibrium queueing-delay pattern and the solution to a linear programming (LP) problem obtained by relaxing some conditions in the original DUE-RDTC problem.
We first present a method for determining whether the GQRP holds.
Based on the GQRP, we then develop a systematic procedure to obtain an exact DUE solution by sequentially solving two LPs: one for the equilibrium cost pattern, including queueing delays, and the other for the corresponding equilibrium flow pattern.
Computational results on networks of varying scales confirm the effectiveness of the proposed method.
}

\def\keyword{dynamic traffic assignment \\
 dynamic user equilibrium \\
 route and departure time choice \\
 linear programming \\ 
 queue replacement principle}
\begin{document}
\maketitle
\titlepageext

\newpage
\section{Introduction}
The dynamic user equilibrium (DUE) assignment with route and departure time choice (RDTC) is a fundamental framework for analyzing dynamic traffic phenomena on networks.
The DUE-RDTC assignment models within-day traffic flow patterns based on the interactions between users' route and departure time choices and the traffic dynamics on road networks.
Clarifying its theoretical properties and developing efficient numerical methods are essential for designing and evaluating traffic management policies.
\par
Studies on the DUE-RDTC assignment have a long history~\citep[see][for detailed reviews]{Han2019-bn}.
Early studies~\citep[e.g.][]{Kuwahara1987-by,Arnott1990-ta} focused on simple network structures, such as Braess networks or networks with a single bottleneck per route.
These studies derived closed-form solutions and examined the properties of equilibrium flow and cost patterns.
However, their analyses relied on the proof-by-cases approach for verifying the equilibrium conditions.
Since the number of cases rapidly increases with the number of bottlenecks, applying the proof-by-cases approach to more general networks is difficult.
\par
Since the 2000s, numerical approaches for general large-scale networks have been developed~\citep{Huang2002-go,Szeto2004-zg,Friesz2006-rs,Friesz2011-jw,Han2019-bn}.
These studies formulated the DUE-RDTC assignment as variational inequality (VI) or differential variational inequality (DVI) problems using the dynamic network loading (DNL) operator~\citep{Friesz2019-ac}.
While their approach established the existence of DUE-RDTC in general networks, it is difficult to derive deeper insights into the theoretical properties of the equilibrium: 
the equilibrium conditions formulated in a conventional Eulerian coordinate system, where state variables are defined with respect to the absolute time, have a complicated nested structure between link and route travel times, and this causes extreme difficulties in analyzing the equilibrium properties~\citep[][]{Akamatsu2015-ip}.
Furthermore, since the cost function is generally non-monotonic, the convergence of fixed-point-based algorithms for VIs and DVIs is typically not guaranteed~\citep{Friesz2011-jw,Friesz2019-ac}.
Indeed, recent numerical studies show that these algorithms often fail to produce accurate equilibrium solutions, particularly on large-scale networks.\footnote{In numerical experiments of \cite{Han2019-bn}, the travel cost gaps for all origin-destination pairs, which should be zero in an exact solution, remained significant in large-scale networks.}
\par
To address the aforementioned difficulties, this study develops a hybrid analytical and numerical method for the DUE-RDTC assignment in many-to-one networks with homogeneous users.
The proposed method is based on the following two key concepts recently developed for the DUE assignment with departure time choice (DTC): formulation in a Lagrangian-like coordinate system and the \textit{queue replacement principle} (QRP).
Before presenting the proposed method, we briefly review these two concepts.
\par
The Lagrangian-like coordinate system approach defines the link-based traffic variables with the departure time at a single destination in a many-to-one network. 
Unlike the conventional Eulerian coordinate system, this definition is well-suited for handling the ex-post trip costs that users actually experience during their trips and for tracing the equilibrium time-space paths of users in a simple and transparent manner~\citep{Kuwahara1993-vt}.
As a result, we can formulate equilibrium problems without the need for complicated nested structures in trip cost functions. 
Moreover, the Lagrangian-like coordinate system provides a significant advantage for the analysis of the DUE assignment when combined with the QRP technique.
\par
The concept of the QRP was formally introduced by \citet{Fu2022-st,Sakai2024-dt,Sakai2024-wy} in the context of the DUE assignment with DTC in the corridor network with multiple bottlenecks.
These studies defined the QRP as a replaceability between the equilibrium queueing delay pattern and the bottleneck-based congestion pricing pattern that eliminates queues on all links.\footnote{
  The concept of the QRP is a well-known property in single bottleneck cases.
  For example, \citet{Vickrey1969-ic} and \citet{Hendrickson1981-cu} presented this property by explaining that congestion pricing mimicking an equilibrium queueing delay pattern eliminates the queue and achieves the dynamic system optimal (DSO) state.
  Conversely, using the QRP, \citet{Iryo2007-ne} and \citet{Akamatsu2021-di} showed that the equilibrium queueing delay pattern can be obtained from the queue-eliminating pricing pattern in the DSO state.
}
This QRP does not always hold; however, when it does, it enables the construction of the DUE state, including the flow and cost patterns, via the queue-eliminating pricing pattern in the DSO state.
The DSO problem is derived by replacing the queueing condition in the original DUE formulation with a simpler bottleneck capacity constraint and is eventually formulated as a linear programming (LP) problem.
The optimal congestion pricing pattern can be obtained as the optimal Lagrange multiplier associated with the bottleneck capacity constraint, i.e., the optimal shadow price.
Thus, this approach mathematically means that the solution to an equilibrium problem can also be derived through the solution to a problem with some of the conditions of the equilibrium problem relaxed.
Although this approach is applicable and powerful, it remains unclear whether an exact DUE solution can be derived for the DUE-RDTC problem, which incorporates route choice, using a similar approach.
\par
This study extends and generalizes the previously proposed \textit{queue replacement approach} for the DUE-DTC model to the DUE-RDTC model.
To achieve this, we introduce the concept of the \textit{Generalized Queue Replacement Principle} (GQRP).
The GQRP establishes an equivalence between the equilibrium queueing-delay pattern and the optimal Lagrange multiplier (shadow price) of an auxiliary LP obtained by relaxing the queueing condition in the original DUE-RDTC formulation.
If the GQRP holds, the equilibrium queueing-delay pattern can be recovered from the optimal Lagrange multiplier, and the corresponding equilibrium flow pattern can then be determined accordingly.\footnote{
The auxiliary problem associated with the DUE-RDTC model is not, in general, a standard DSO problem, except in special cases, unlike the DUE-DTC model considered in previous studies.
Nevertheless, our analysis shows that, under the GQRP, the same shadow-price mechanism still characterizes the equilibrium queueing pattern.
}
We develop a systematic methodology for verifying whether the GQRP holds and constructing an exact DUE solution including equilibrium flow and cost patterns based on this GQRP.
\begin{figure}[tbp]
  \centering
  \includegraphics[width=0.8\textwidth]{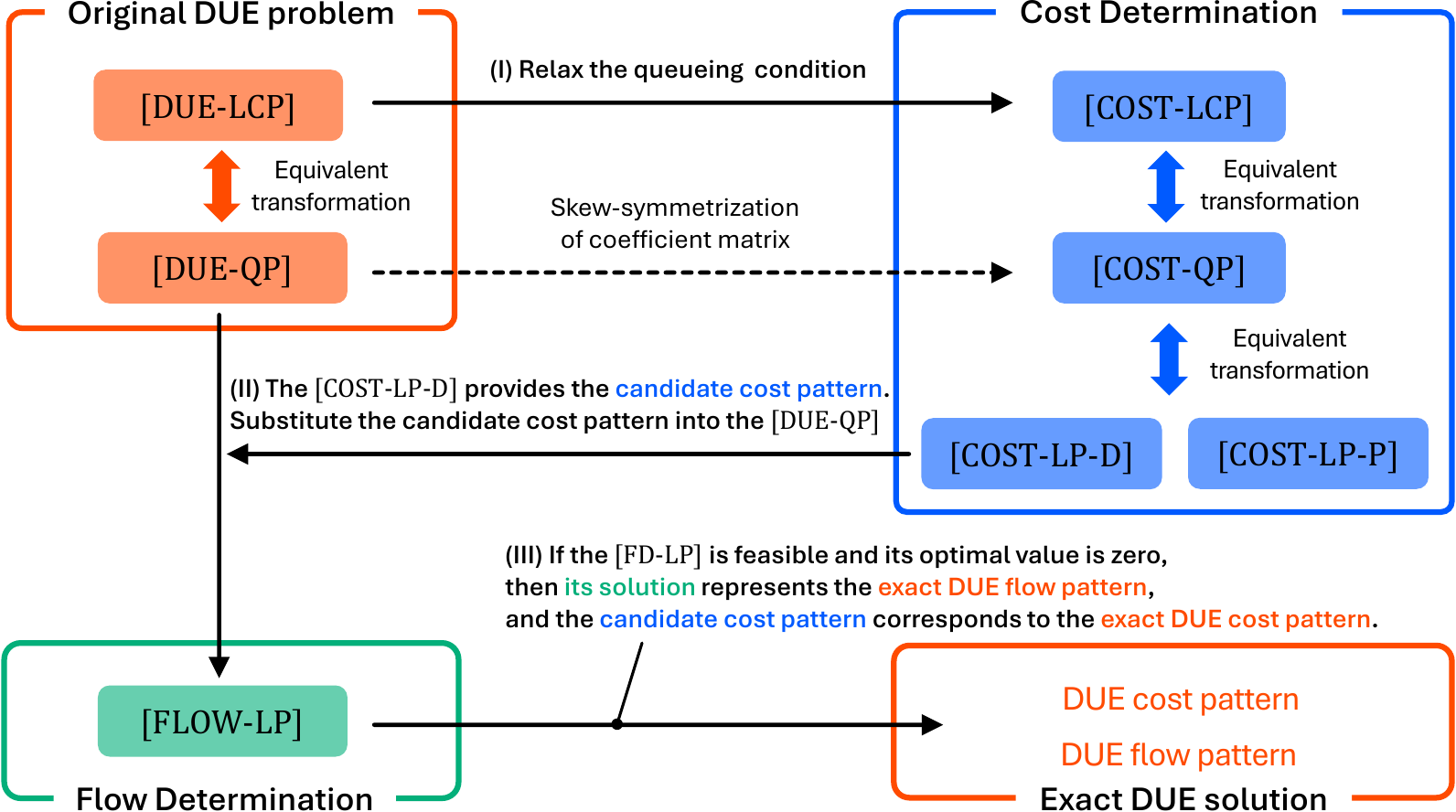}
  \caption{Framework of the proposed GQRP-based approach}
  \label{fig:overview}
\end{figure}
\par
Technically, the proposed approach consists of two main steps: 
(i) the cost determination step and 
(ii) the GQRP verification and flow determination step (see \cref{fig:overview}). 
The approach begins by formulating the DUE-RDTC problem as a mixed linear complementarity problem (LCP), referred to as the \texttt{[DUE-LCP]}, using the Lagrangian-like coordinate system. 
In this formulation, the unknowns are the equilibrium travel cost pattern (including queueing delays) and the equilibrium flow pattern.
The \texttt{[DUE-LCP]} can be transformed into an equivalent quadratic program (\texttt{[DUE-QP]}), but this program has an asymmetric structure that is difficult to solve directly.
To address this, in the cost determination step, we relax the queueing condition\footnote{
  Strictly speaking, we also relax the consistency condition, which is a condition that guarantees the Lipschitz continuity of cumulative arrival flows.
} and obtain a skew-symmetric corresponding QP problem, referred to as the \texttt{[COST-QP]}.
This skew-symmetry allows us to decompose the problem into two linear programs: a primal LP (\texttt{[COST-LP-P]}) and a dual LP (\texttt{[COST-LP-D]}).
Under certain conditions, the solution of \texttt{[COST-LP-D]} gives the exact equilibrium travel cost pattern.
This property is the essence of the GQRP.
\par
In the second step, we verify whether the GQRP holds (i.e., whether the candidate cost pattern obtained from \texttt{[COST-LP-D]} is an exact equilibrium travel cost pattern) and determine the equilibrium flow pattern.
Specifically, we substitute the candidate cost pattern into the original DUE formulation.
This substitution reduces the problem to a flow determination LP (\texttt{[FLOW-LP]}).
We find that if the \texttt{[FLOW-LP]} is feasible and its optimal value is zero, then the candidate cost pattern corresponds to the exact equilibrium travel cost pattern, and the solution of the \texttt{[FLOW-LP]} gives the exact equilibrium flow pattern.
Thus, the DUE-RDTC solution can be obtained by sequentially solving two LPs.
Both LPs are computationally tractable, providing a practical and efficient alternative to directly solving the original complex formulation.
{
\par
Compared with the previous QRP-based studies \citep{Fu2022-st, Sakai2024-dt, Sakai2024-wy}, the present study extends the queue replacement approach in two respects.
First, regarding the model setting, the previous studies addressed only the DUE-DTC model on a corridor network with a single route, whereas the present study establishes the GQRP for the DUE-RDTC model on many-to-one networks with simultaneous route and departure time choice.
Second, regarding the underlying principle, the previous studies relied on the equivalence between the auxiliary LP and a standard DSO problem, whereas the present study generalizes the QRP into a principle that does not require such an equivalence: although the auxiliary LP for the DUE-RDTC model is not, in general, a standard DSO problem, the same shadow-price (Lagrange multiplier) mechanism still characterizes the equilibrium queueing-delay pattern under the GQRP.
}
\par
The paper is organized as follows.
\cref{sec:Model} formulates the DUE-RDTC problem as a mixed LCP using a Lagrangian-like coordinate system.
\cref{sec:QueueReplacement} introduces the queue replacement approach based on the GQRP.
\cref{sec:Discussion} examines the policy and algorithmic implications of the proposed framework.
\cref{sec:IllustrativeExamples} presents numerical experiments demonstrating its effectiveness.
Finally, \cref{sec:Conclusion} summarizes the findings and concludes the paper.

\section{Model}
\label{sec:Model}
\subsection{Network and user}
\label{sec:NetworkAndUser}
\begin{figure}[tbp]
  \centering
  \begin{minipage}{0.50\textwidth}
    \centering
    \includegraphics[width=0.7\textwidth]{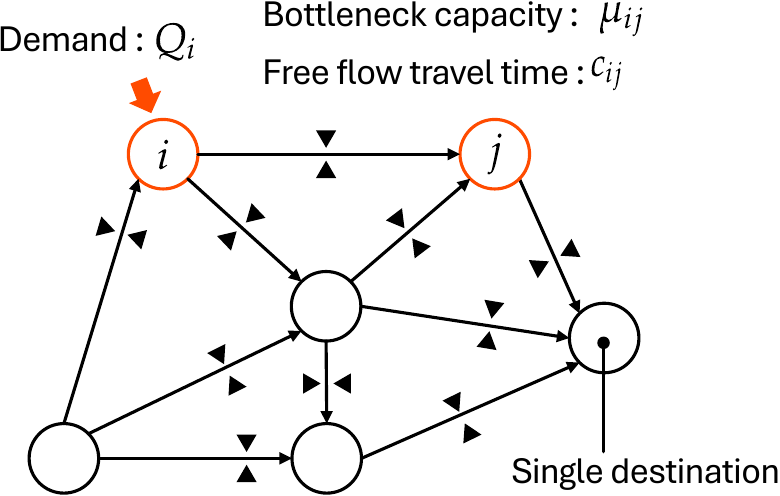}
    \caption{Many-to-one network}
    \label{fig:network}
  \end{minipage}
  \centering
  \begin{minipage}{0.45\textwidth}
    \centering
    \includegraphics[width=0.8\textwidth]{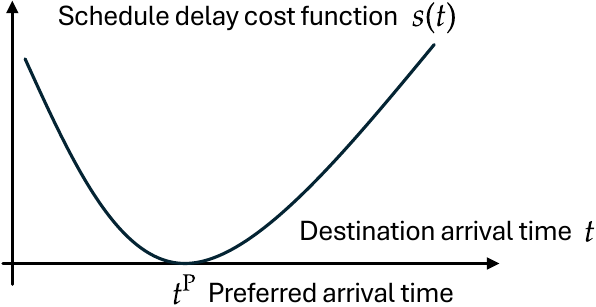}
    \caption{Schedule delay cost function}
    \label{fig:sdcf}
  \end{minipage}
\end{figure}

Consider a network with $N+1$ nodes and $L$ directed links as shown in \cref{fig:network}.
The nodes are classified into two types: origin nodes and a single destination node.
The set of origin nodes and links are denoted by $\ClN$ and $\ClL$, respectively.
All links consist of the bottleneck section and the free flow section.
The bottleneck capacity and the free flow travel time of link $(i, j) \in \ClL$ are denoted by $\mu_{ij}$ and $c_{ij}$, respectively.
A queue is formed at each bottleneck when the inflow rate exceeds the bottleneck capacity. 
The queue dynamics are modeled by the standard point queue model along with the first-in-first-out (FIFO) principle.
\par
The demand from origin $i$ to the destination is denoted by $Q_{i}$.
Each user incurs the free flow travel cost and the queueing delay cost on each link on their route.
Additionally, the user incurs the schedule delay cost at the destination, which stems from the difference between the preferred arrival time and the actual arrival time.
We assume that all users have the same preferred arrival time at the destination, which is denoted by $t^{\mathrm{P}}$.
Users who arrive at the destination at time $t$ incur a cost of $s(t \mid t^{\mathrm{P}})$.
In the following, $s(t \mid t^{\mathrm{P}})$ is written as $s(t)$ for simplicity and called the schedule delay cost function.
The schedule delay cost function $s(t)$ is assumed to be convex in $t$ and $s(t^{\mathrm{P}}) = 0$, as shown in \cref{fig:sdcf}.
Following previous studies, we also assume $\mathrm{d} s(t) / \mathrm{d} t > -1$ for the validity of the departure time choice model~\citep{Daganzo1985-ls,Lindsey2004-aw}.
In addition, we assume that $\mathrm{d}s(t)/\mathrm{d}t$ is finite for all $t$.

\subsection{Dynamic user equilibrium}
\begin{figure}
    \centering
    \includegraphics[width=0.7\linewidth]{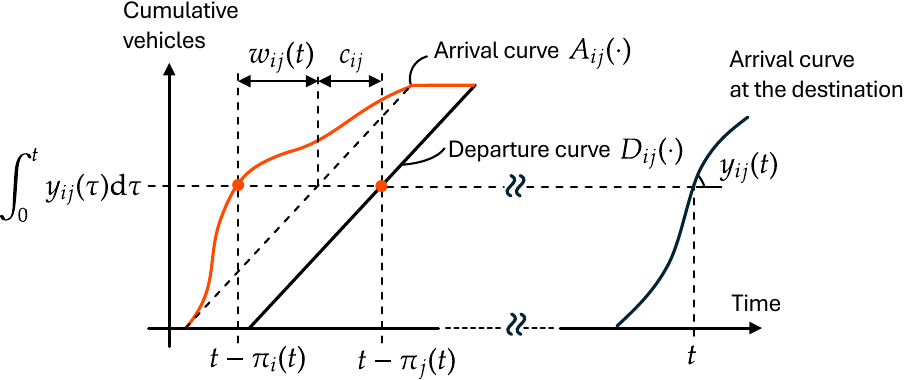}
    \caption{Definition of variables in the Lagrangian-like coordinate system}
    \label{fig:LtoE}
\end{figure}

This section formulates the DUE problem by using the \textit{Lagrangian-like coordinate system} approach.
We formally introduce the definition of the DUE state as follows.
\begin{dfn}
  \label{dfn:DUE}
  The DUE state is a state in which each user chooses a departure time and route such that their experienced travel cost is minimal, and no alternative route or departure time provides a lower cost, given congestion caused by all other users.
\end{dfn}
\noindent
From \cref{dfn:DUE}, all users have no incentive to change their routes or departure times in the DUE state.
This means that, in the DUE state, users who arrive at the destination at the same time have the same departure times at the nodes on their routes~\citep{Friesz1993-bj,Kuwahara1993-vt,Akamatsu2001-bi}.
Thus, users' arrival times at the pass-through nodes are uniquely determined by their destination arrival time.
Based on the property, in a Lagrangian-like coordinate system approach, the state variables are defined in association with the arrival time at the destination~(see \citet{Kuwahara1993-vt,Akamatsu2000-ph,Akamatsu2015-ip,Akamatsu2021-di} for the detailed discussion).
Using this approach, we define the state of the variables as follows:
\begin{description}
    \item[] $q_{i}(t)$ : destination arrival flow of users whose origin is $i$ and whose destination arrival time is $t$ 
    \item[] $y_{ij}(t)$ : destination arrival flow of users who pass through link $(i, j)$ and {arrive} at the destination at time $t$ 
    \item[] $w_{ij}(t)$ : queueing delay on link $(i, j)$ for users who arrive at the destination at time $t$ 
    \item[] $\pi_{i}(t)$ : earliest travel time from node $i$ to the destination for a user who {arrives} at the destination at time $t$ 
\end{description}
Note that $y_{ij}(t)$ is the flow in the Lagrangian coordinate system, i.e., it represents the arrival flow at the destination at time $t$.
To obtain the flow in the Eulerian coordinate system from $y_{ij}(t)$, we trace the time-space trajectories of users by using $\pi_{i}(t)$ and $\pi_{j}(t)$ as shown in \cref{fig:LtoE}.
The figure represents the relationships between the variables ($y_{ij}(t)$, $w_{ij}(t)$, $\pi_{i}(t)$) and the variables in the Eulerian coordinate system (i.e., cumulative arrival curve $A_{ij}(\cdot)$ and departure curve $D_{ij}(\cdot)$) at the link $(i, j)$.
The figure shows that we obtain the cumulative arrival curve $A_{ij}(\cdot)$ and departure curve $D_{ij}(\cdot)$ based on the arrival flow at the destination $y_{ij}(t)$ and the travel times  $\pi_{i}(t)$ and $\pi_{j}(t)$.
\par
\noindent
In the Lagrangian-like coordinate system, the travel cost of a user who arrives at the destination at time $t \in \ClT$ with route $\ClP$ is represented as follows:
\begin{align}
  C(t \mid \ClP) = \sum_{(i, j) \in \ClP} \alpha \left( c_{ij} + w_{ij}(t) \right) + s(t),
\end{align}
where $\alpha$ is the parameter of the value of time and is assumed to be $1$ for simplicity.
\par
Using these variables, the DUE problem can be formulated as a mixed linear complementarity problem (LCP) consisting of six conditions.
Hereinafter, we use matrix and vector notation, represented by bold symbols (e.g., $\Vtq(t)$, $\Vtmu$).
\par
\noindent
\textbf{(i) Demand conservation condition:}
This condition ensures that the total commuting flow from each origin is equal to the demand from each origin as follows:
\begin{align}
  &\int_{t \in \ClT} \Vtq(t) \mathrm{d} t= \VtQ.
  \label{eq:demand_conservation}
\end{align}
The element-wise representation of \cref{eq:demand_conservation} is as follows:
\begin{align}
  &\int_{t \in \ClT} q_{i}(t) \mathrm{d}t = Q_{i}
  &&\forall i \in \ClN.
  \tag{2a}
\end{align}
\par
\noindent
\textbf{(ii) Flow conservation condition:}
This condition ensures that, at each node, the inflow and the demand flow are equal to the outflow.
\begin{align}
  &\VtA \Vty(t) - \Vtq(t) = \Vt0
  &&\forall t \in \ClT,
  \label{eq:flow_conservation}
\end{align}
where $\VtA \in \mathbb{R}^{N \times L}$ is the node-link incidence matrix excluding the destination node.
The element of the $k$-th row and $ij$-th column of $\VtA$ is $1$ if node $k$ is the upstream end of link $(i, j)$ (i.e., link $(i, j)$ flows out from node $k$), $-1$ if node $k$ is the downstream end of link $(i, j)$ (i.e., link $(i, j)$ flows into node $k$), and $0$ otherwise. 
The element-wise representation of \cref{eq:flow_conservation} is as follows:
\begin{align}
  &\sum_{j \in \ClN^{\mathrm{OUT}}(i)} y_{ij}(t) = 
  \sum_{j \in \ClN^{\mathrm{IN}}(i)} y_{ji}(t) + q_{i}(t)
  &&\forall i \in \ClN, \quad \forall t \in \ClT,
  \tag{3a}
\end{align}
where $\ClN^{\mathrm{OUT}}(i)$ is denoted by the set of nodes that flow out from node $i$ and $\ClN^{\mathrm{IN}}(i)$ is denoted by the set of nodes that flow into node $i$.
\par
\noindent
\textbf{(iii) Queuing condition:}
This condition represents the point queue dynamic at the bottleneck.
The condition is expressed as the following linear complementarity condition~\citep{Akamatsu2000-ph,Akamatsu2015-ip}:
\begin{align}
  &\Vt0 \leq 
  \left(
      \VtM \dot{\Vtw}(t) - \VtM \VtA^{\top}_{+} \dot{\Vtpi}(t) + \Vtmu
   - \Vty(t)
  \right) \perp \Vtw(t) \geq \Vt0
  &&\forall t \in \ClT,
  \label{eq:queueing}
\end{align}
where $\VtM$ is a diagonal matrix whose diagonal elements are $\{ \mu_{ij} \}$ and $\VtA_{+}$ is the node-link incidence matrix $\VtA$ with all negative elements replaced by $0$.
The overdot denotes the derivative of the variable with respect to the destination arrival time $t \in \ClT$.
The element-wise representation of the queueing condition is as follows:
\begin{align}
  &\begin{dcases}
    \mu_{ij} \left( \dot{w}_{i j}(t)
    -\dot{\pi}_{i}(t) + 1 \right)
    - y_{ij}(t)
    = 0
    &\text{if} \quad w_{ij}(t) > 0
    \\
    \mu_{ij} \left( \dot{w}_{i j}(t)
    -\dot{\pi}_{i}(t) + 1 \right)
    - y_{ij}(t)
    \geq 0
    &\text{if} \quad w_{ij}(t) = 0
  \end{dcases}
  &&\forall (i, j) \in \ClL, \quad \forall t \in \ClT.
    \tag{4a}
\end{align}
This condition represents that if a bottleneck queue exists, then the bottleneck outflow is equal to the capacity $\mu_{ij}$ (see \cref{fig:LtoE}).
In the usual representation without the Lagrangian coordinate system, this queueing condition means $\mathrm{d} D_{ij}(u) / \mathrm{d} u = \mu_{ij}$ if queue exists at the bottleneck, where $D_{ij}(u)$ is the cumulative outflow curve at link $(i, j)$.
\par
\noindent
\textbf{(iv) Consistency condition:}
This condition guarantees the Lipschitz continuity of cumulative arrival flows, which means that the traffic flow is physically consistent~\citep{Akamatsu2015-ip,Fu2022-st,Sakai2022-vm}.
This condition is expressed as follows:
\begin{align}
  &\dot{\Vtpi}(t) \leq \Vt1
  &&\forall t \in \ClT.
  \label{eq:consistency}
\end{align}
The element-wise representation of \cref{eq:consistency} is $\dot{\pi}_{i}(t) \leq 1$,  $\forall i \in \ClN, \quad \forall t \in \ClT$.
\par
\noindent
\textbf{(v) Route choice condition:}
In the DUE state, all users have no incentive to change their routes.
This condition is expressed by the following linear complementarity condition~\citep{Akamatsu2000-ph}:
\begin{align}
 &\Vt0 \leq 
  \left(
    \Vtw(t) - \VtA^{\top} \Vtpi(t) + \Vtc
  \right) \perp \Vty(t) \geq \Vt0
 &&\forall t \in \ClT.
 \label{eq:route_choice}
\end{align}
The vector $\Vtc$ is the vector representation of the free flow travel time $c_{ij}$ for each link $(i, j) \in \ClL$.
The element-wise representation of \cref{eq:route_choice} is as follows:
\begin{align}
  &\begin{dcases}
    w_{ij}(t) - \pi_{i}(t) + \pi_{j}(t) + c_{ij} = 0
    &\text{if} \quad y_{ij}(t) > 0
    \\
    w_{ij}(t) - \pi_{i}(t) + \pi_{j}(t) + c_{ij} \geq 0
    &\text{if} \quad y_{ij}(t) = 0
  \end{dcases}
  &&\forall (i, j) \in \ClL, \quad \forall t \in \ClT.
  \tag{6a}
\end{align}
\par
\noindent
\textbf{(vi) Departure time choice condition:}
In the DUE state, all users have no incentive to change their departure times.
That is, all users from the origin node $i$ incur the same commuting cost.
This condition is expressed as follows:
\begin{align}
  &\Vt0 \leq \left(
    \Vtpi(t) + \Vts(t) - 
    \Vtrho
  \right) \perp \Vtq(t) \geq \Vt0
  &&\forall t \in \ClT,
  \label{eq:departure_time_choice}
\end{align}
where $\Vtrho$ is the vector representation of $\rho_{i}$,  which is the equilibrium travel cost of the origin $i$ user.
The element-wise representation of \cref{eq:departure_time_choice} is as follows:
\begin{align}
  &\begin{dcases}
    \pi_{i}(t) + s(t) - \rho_{i} = 0
    &\text{if} \quad q_{i}(t) > 0
    \\
    \pi_{i}(t) + s(t) - \rho_{i} \geq 0
    &\text{if} \quad q_{i}(t) = 0
  \end{dcases}
  &&\forall i \in \ClN, \quad \forall t \in \ClT.
  \tag{7a}
\end{align}
\par
By combining these conditions, the DUE problem can be formulated as a mixed linear complementarity problem (LCP) as follows:
\begin{align}
\texttt{[DUE-LCP]} \quad \quad \quad \quad
\text{Find}
\quad
&\Vtx \equiv
\left[
  \Vtq,
  \Vty,
  \Vtw,
  \Vtpi,
  \Vtrho
\right] \geq \Vt0,
\quad
\text{such that} \quad
\text{\cref{eq:demand_conservation,eq:flow_conservation,eq:route_choice,eq:departure_time_choice,eq:queueing,eq:consistency}.}
\end{align}
Here, we define 
$\Vtq \equiv \left\{ \Vtq(t) \right\}_{t \in \ClT}$, 
with $\Vty$, $\Vtw$, and $\Vtpi$ being defined in the same manner.
The solution to the \texttt{[DUE-LCP]} is represented by $\UE{\Vtx}$.
Note that, at some times, the value of $\pi_{i}(t)$ is not uniquely determined by the equilibrium conditions alone and may therefore be indeterminate~\citep{Akamatsu2000-ph}.
In such cases, $\pi_{i}(t)$ can be specified consistently, in accordance with its definition, by computing the earliest travel time implied by the given queueing delay pattern $\Vtw$.
\par
The \texttt{[DUE-LCP]} can be transformed into a quadratic programming (QP) problem by using the following lemma, which is a general property of LCPs.
\begin{lem}[\citet{Cottle2009-xz}]
  \label{lem:LCPtoQP}
  Consider a following LCP and QP:
  \begin{align}
  &\texttt{[LCP]} \quad \text{Find} \quad \Vtx^{\ast} \quad \text{such that} \quad \Vt0 \leq \VtB\Vtx^{\ast} + \Vtb \perp \Vtx^{\ast} \geq \Vt0.
  \\
  &\texttt{[QP]} \quad \min_{\Vtx}. \quad f(\Vtx) = \Vtx^{\top}\VtB\Vtx + \Vtb^{\top}\Vtx  \quad \text{s.t.} \quad \VtB\Vtx + \Vtb \geq \Vt0, \quad \Vtx \geq \Vt0.
  \end{align}
  Let $\Vtx^{\ast\ast}$ be a solution to the \texttt{[QP]}.
  The solution to the \texttt{[LCP]} is given by $\Vtx^{\ast} = \Vtx^{\ast\ast}$ and the optimal value of the \texttt{[QP]} is zero, i.e., $f(\Vtx^{\ast\ast})=0$.
\end{lem}
By using \cref{lem:LCPtoQP}, we transform the \texttt{[DUE-LCP]} into a QP with the consistency condition.
This is shown in the following proposition.
\begin{pro}
\label{pro:DUELCPtoDUEQP}
The solution to the \texttt{[DUE-LCP]} is given by the solution to the following QP problem:
\begin{align}
\texttt{[DUE-QP]} \quad \min_{\left[ \Vtq, \Vty, \Vtw, \Vtpi, \Vtrho \right] \geq \Vt0}. \quad
&\UE{Z}(\Vtq, \Vty, \Vtw, \Vtpi, \Vtrho)
\equiv Z^{\mathrm{Q}}(\Vtw, \Vtpi) 
       + Z^{\mathrm{L}}(\Vtq, \Vty, \Vtw, \Vtpi, \Vtrho)
\\
\text{s.t.} \quad
&\int_{t \in \ClT} \Vtq(t) \mathrm{d}t - \VtQ = \Vt0,
\\
&\VtA \Vty(t) - \Vtq(t) = \Vt0
&&\forall t \in \ClT,
\\
& \left(
  \VtM \dot{\Vtw}(t) - \VtM \VtA^{\top}_{+} \dot{\Vtpi}(t) + \Vtmu
- \Vty(t)
\right) \geq \Vt0
&&\forall t \in \ClT,
\label{eq:DUE-QP_queueing}
\\
&\dot{\Vtpi}(t) \leq \Vt1
&&\forall t \in \ClT,
\\
&\Vtw(t) - \VtA^{\top} \Vtpi(t) + \Vtc \geq \Vt0
&&\forall t \in \ClT,
\\
& \Vtpi(t) + \Vts(t) - \Vtrho \geq \Vt0
&&\forall t \in \ClT,
\end{align}
\begin{align}
\text{where} \quad
&Z^{\mathrm{Q}}(\Vtw, \Vtpi)
\equiv \int_{t \in \ClT} \Vtw(t) \left( \VtM \dot{\Vtw}(t) \right)\mathrm{d} t
- \int_{t \in \ClT} \Vtw(t) \left( \VtM \VtA^{\top}_{+} \dot{\Vtpi}(t) \right)\mathrm{d} t,
\\
&Z^{\mathrm{L}}(\Vtq, \Vty, \Vtw, \Vtpi, \Vtrho)
\equiv  
Z^{\mathrm{P}}(\Vtq, \Vty) -  Z^{\mathrm{D}}(\Vtw, \Vtpi, \Vtrho)
\\
&Z^{\mathrm{P}}(\Vtq, \Vty)
\equiv \int_{t \in \ClT} \Vts(t) \Vtq(t) \mathrm{d} t
+ \int_{t \in \ClT} \Vtc \Vty(t) \mathrm{d}t,
\\
&Z^{\mathrm{D}}(\Vtw, \Vtpi, \Vtrho) \equiv
- \int_{t \in \ClT} \Vtmu \Vtw(t) \mathrm{d}t
+ \VtQ^{\top} \Vtrho.
\end{align}
The optimal value of the \texttt{[DUE-QP]} is zero, i.e., $\UE{Z}(\UE{\Vtq}, \UE{\Vty}, \UE{\Vtw}, \UE{\Vtpi}, \UE{\Vtrho}) = 0$.
\end{pro}
The objective function $\UE{Z}(\Vtq, \Vty, \Vtw, \Vtpi, \Vtrho)$ consists of a quadratic term $Z^{\mathrm{Q}}(\Vtw, \Vtpi)$ and a linear term $Z^{\mathrm{L}}(\Vtq, \Vty, \Vtw, \Vtpi, \Vtrho)$.
The quadratic term $Z^{\mathrm{Q}}(\Vtw, \Vtpi)$ includes differential and quadratic expressions in $\Vtw$ and $\Vtpi$.
The linear term $Z^{\mathrm{L}}(\Vtq, \Vty, \Vtw, \Vtpi, \Vtrho)$ can be further decomposed into $Z^{\mathrm{P}}(\Vtq, \Vty)$ and $Z^{\mathrm{D}}(\Vtw, \Vtpi, \Vtrho)$, which depend on the flow variables $[\Vtq(t), \Vty(t)]$ and the cost variables $[\Vtw, \Vtpi, \Vtrho]$, respectively.
The presence of $Z^{\mathrm{Q}}(\Vtw, \Vtpi)$ makes the problem non-convex, preventing the use of standard convex optimization techniques.
Moreover, the constraint~\cref{eq:DUE-QP_queueing} contains derivatives of $\Vtw$ and $\Vtpi$, introducing asymmetry and making the direct solution for the \texttt{[DUE-QP]} more complicated.

\section{Queue Replacement Approach}
\label{sec:QueueReplacement}
\subsection{Overview of the approach}
This section presents the queue replacement approach to the DUE-RDTC problem.
The core idea is the Generalized Queue Replacement Principle (GQRP), which states that, under certain conditions, the equilibrium queueing delay patterns in the original DUE problem can be obtained by solving an auxiliary LP problem derived by relaxing certain constraints of the original problem.
Because LPs can be solved efficiently using standard numerical techniques, the GQRP can significantly reduce the computational complexity typically associated with solving the DUE-RDTC problem, which inherently involves complex LCPs.
The proposed approach consists of two main sequential steps (see \cref{fig:overview} for a graphical overview):

\paragraph{Step 1: Cost Determination}
In this step, we derive a candidate equilibrium queueing delay pattern.
The cost determination problem for obtaining this pattern is formulated by relaxing the consistency condition and replacing the original DUE problem's queueing condition with a simpler one.
This relaxation yields a QP with a skew-symmetric structure, which can be transformed into a primal-based LP \texttt{[COST-LP-P]} and a dual-based LP \texttt{[COST-LP-D]}.
The solution to the \texttt{[COST-LP-D]} then provides the candidate queueing delay pattern.

\paragraph{Step 2: Flow Determination and the GQRP Verification}
In this step, we derive the equilibrium flow pattern and verify whether the candidate queueing delay pattern obtained in the first step satisfies the original DUE conditions.  
To do so, we substitute the candidate pattern into the original DUE formulation, i.e., into \texttt{[DUE-QP]}.  
This substitution yields a reduced flow-determination problem, formulated as a linear program \texttt{[FLOW-LP]}.  
If \texttt{[FLOW-LP]} is feasible and its optimal objective value is zero, the candidate queueing delay pattern is indeed the exact equilibrium queueing delay pattern.  
Moreover, the solution to \texttt{[FLOW-LP]} provides the exact equilibrium flow pattern.

\subsection{Cost determination step}
\label{sec:CostDetermination}
The cost determination problem is derived from the original DUE problem by removing the consistency condition \eqref{eq:consistency} and relaxing the queueing condition \eqref{eq:queueing}.
Specifically, we exchange the queueing condition \eqref{eq:queueing} with the following condition:
\begin{align}
  &\begin{dcases}
    \mu_{ij}
    - y_{ij}(t)
    = 0
    &\text{if} \quad w_{ij}(t) > 0
    \\
    \mu_{ij} 
    - y_{ij}(t)
    \geq 0
    &\text{if} \quad w_{ij}(t) = 0
  \end{dcases}
  \hspace{2cm} \forall (i, j) \in \ClL, \quad \forall t \in \ClT.
\end{align}
The vector representation of the above condition is as follows:
\begin{align}
  &\Vt0 \leq \Vtmu - \Vty(t) \perp \Vtw(t) \geq \Vt0
  &&\forall t \in \ClT.
  \label{eq:CD_queueing_condition}
\end{align}
This representation allows us to avoid dealing with differential variables, i.e., $\dot{\Vtsigma}(t)$ and $\dot{\Vtw}(t)$, which take temporally asymmetric structure.
As a result, we can express the problem as a LCP without differential variables as follows:
\begin{align}
\text{\texttt{[COST-LCP]}} \quad
\text{Find}
\quad
\Vtx \equiv
\left[
  \Vtq(t),
  \Vty(t),
  \Vtw(t),
  \Vtpi(t),
  \Vtrho
\right]^{\top} \geq \Vt0,
\quad
\text{such that} \quad
\text{\cref{eq:demand_conservation,eq:flow_conservation,eq:route_choice,eq:departure_time_choice,eq:CD_queueing_condition}.}
\end{align}
By applying the same transformation as in \cref{pro:DUELCPtoDUEQP} based on \cref{lem:LCPtoQP}, we can transform the \texttt{[COST-LCP]} into a QP problem, referred to as the \texttt{[COST-QP]}.
\begin{pro}
\label{pro:CDLCPtoCDQP}
The solution to the \texttt{[COST-LCP]} is given by the solution to the following QP problem:
\begin{align}
\texttt{[COST-QP]} \quad \min_{\left[ \Vtq, \Vty, \Vtw, \Vtpi, \Vtrho \right] \geq \Vt0}. \quad
&\CD{Z}(\Vtq, \Vty, \Vtw, \Vtpi, \Vtrho) \equiv 
Z^{\mathrm{P}}(\Vtq, \Vty) - Z^{\mathrm{D}}(\Vtw, \Vtpi, \Vtrho) 
\\
\text{s.t.} \quad
&\int_{t \in \ClT} \Vtq(t) \mathrm{d}t - \VtQ = \Vt0,
\\
&\VtA \Vty(t) - \Vtq(t) = \Vt0
&&\forall t \in \ClT,
\\
& \Vtmu - \Vty(t)
\geq \Vt0
&&\forall t \in \ClT,
\\
& \Vtpi(t) + \Vts(t) - \Vtrho \geq \Vt0
&&\forall t \in \ClT,
\\
&  \Vtw(t) - \VtA^{\top} \Vtpi(t) + \Vtc \geq \Vt0
&&\forall t \in \ClT.
\end{align}
\end{pro}
Unlike the \texttt{[DUE-QP]}, the objective function and constraints of the \texttt{[COST-QP]} consist solely of linear terms and contain no differential or quadratic terms.  
This simplification enables the \texttt{[COST-QP]} to be transformed into a primal-based LP (\texttt{[COST-LP-P]}) and a dual-based LP (\texttt{[COST-LP-D]}), as summarized in the following proposition.
\begin{pro}
\label{pro:CDQPtoCDLP}
The solution to the \texttt{[COST-QP]} is given by the solution to the following primal and dual LP problems:
\begin{align}
  \text{\texttt{[COST-LP-P]}} \quad 
  \min_{[\Vtq, \Vty] \geq \Vt0}. \quad & 
  Z^{\mathrm{P}}(\Vtq, \Vty) 
  \\
  \text{s.t.} \quad 
  &\int_{t \in \ClT} \Vtq(t) \mathrm{d}t - \VtQ = \Vt0,
  \label{eq:DSO_demand_conservation}
  \\
  &\VtA \Vty(t) - \Vtq(t) = \Vt0
  &&\forall t \in \ClT,
  \label{eq:DSO_flow_conservation}
  \\
  & \Vty(t) - \Vtmu \leq \Vt0
  &&\forall t \in \ClT.
  \label{eq:DSO_non_queueing}
\end{align}
\begin{align}
  \text{\texttt{[COST-LP-D]}} \quad 
  \max_{[\Vtw, \Vtpi, \Vtrho] \geq \Vt0}. \quad &
  Z^{\mathrm{D}}(\Vtw, \Vtpi, \Vtrho)
  \\
  \text{s.t.} \quad
  & \Vtpi(t) + \Vts(t) - \Vtrho \geq \Vt0
  &&\forall t \in \ClT,
  \\
  &  \Vtw(t) - \VtA^{\top} \Vtpi(t) + \Vtc \geq \Vt0
  &&\forall t \in \ClT.
\end{align}
\end{pro}
This proposition implies that the solution to the \texttt{[COST-LCP]} can be obtained by solving the \texttt{[COST-LP-P]} and the \texttt{[COST-LP-D]}.  
The equivalence among the \texttt{[COST-LCP]}, \texttt{[COST-LP-P]}, and \texttt{[COST-LP-D]} can also be verified by comparing the conditions of the \texttt{[COST-LCP]} with the optimality conditions of the \texttt{[COST-LP-P]} and \texttt{[COST-LP-D]}.  
These optimality conditions are as follows:
\begin{align}
  &\int_{t \in \ClT} \CD{\Vtq}(t) \mathrm{d} t = \VtQ,
  \label{eq:CD_demand_conservation}
  \\
  &\VtA \CD{\Vty}(t) - \CD{\Vtq}(t) = \Vt0
  &&\forall t \in \ClT,
  \label{eq:CD_flow_conservation}
  \\
  &\Vt0 \leq 
  \left(
  \Vtmu - \CD{\Vty}(t)
  \right) \perp \CD{\Vtw}(t) \geq \Vt0
  &&\forall t \in \ClT,
  \label{eq:CD_queue}
  \\
  &\Vt0 \leq \left(
    \CD{\Vtpi}(t) + \Vts(t) - \CD{\Vtrho}
  \right) \perp \CD{\Vtq}(t) \geq \Vt0
  &&\forall t \in \ClT,
  \label{eq:CD_departure}
  \\
  &\Vt0 \leq 
  \left(
    \CD{\Vtw}(t) - \VtA^{\top} \CD{\Vtpi}(t) + \Vtc
  \right) \perp \CD{\Vty}(t) \geq \Vt0
  &&\forall t \in \ClT,
  \label{eq:CD_route}
\end{align}
where $\CD{\Vtx}$ represents the solution to the \texttt{[COST-LP-P]} and \texttt{[COST-LP-D]}.
The solution $\CD{\Vtw}(t)$ obtained from the \texttt{[COST-LP-D]} is called a \textit{candidate} queueing delay pattern because it does not necessarily satisfy the original queueing condition~\eqref{eq:queueing}.
Whether this $\CD{\Vtw}(t)$ is the true equilibrium queueing delay pattern will be confirmed in the subsequent flow determination step.
\par
Here, we remark that the \texttt{[COST-LP-P]} also plays an important role in the proposed approach.
First, the \texttt{[COST-LP-P]} has the structure of a minimum-cost flow problem with capacity constraints that are independent of queueing delay.\footnote{
The relationship between the \texttt{[COST-LP-P]} and the dynamic system optimal problem is discussed in \cref{sec:FirstBestPolicies}.
}
Using this structure, primal-dual-based algorithms that utilize the \texttt{[COST-LP-P]} may offer computational improvements over dual-based methods that rely solely on the \texttt{[COST-LP-D]}.
Second, in certain special cases, part of the solution to the \texttt{[COST-LP-P]} coincides with part of the exact DUE flow pattern.
This aspect will be discussed in \cref{sec:ConstructiveApproach} in the context of the subsequent flow determination step.

\subsection{Flow determination step}
\label{sec:FlowDetermination}

\subsubsection{Generalized queue replacement principle and numerical flow determination}
We test whether the GQRP holds—i.e., whether the candidate queueing delay pattern $\CD{\Vtw}(t)$ from the first step is an exact equilibrium pattern—by checking for the existence of a flow pattern that satisfies the original DUE conditions.
The discussion begins by formal definition of the GQRP as follows:
\begin{dfn}[Generalized Queue Replacement Principle]
  \label{dfn:QRP}
  If the queueing delay pattern $\UE{\Vtw}(t)$ in the \texttt{[DUE-LCP]} is equal to the solution $\CD{\Vtw}(t)$ in the \texttt{[COST-LCP]}, that is,
  \begin{align}
    &\UE{\Vtw}(t) = \CD{\Vtw}(t)
    &&\forall t \in \ClT,
    \label{eq:cnj_w}
  \end{align}
  then the generalized queue replacement principle holds.
\end{dfn}
The flow determination problem is obtained by substituting the candidate queueing delay pattern $\CD{\Vtw}(t)$ into the original DUE formulation \texttt{[DUE-LCP]} (effectively, \texttt{[DUE-QP]}).
In the substitution process, it becomes clear that $\UE{\Vtpi}(t) = \CD{\Vtpi}(t)$, since the DUE route choice condition~\eqref{eq:route_choice} shares the same structure as the optimality condition~\eqref{eq:CD_route}.
Similarly, from the departure time choice condition~\eqref{eq:departure_time_choice} and the corresponding optimality condition~\eqref{eq:CD_departure}, we see $\UE{\Vtrho} = \CD{\Vtrho}$.
Furthermore, regarding the consistency condition~\eqref{eq:consistency} in the original DUE problem, we obtain the following lemma:
\begin{lem}
  \label{lem:CD_pi_dot<1}
  The solution $\CD{\Vtpi}(t)$ satisfies the following condition:
  \begin{align}
    &\CD{\dot{\pi}}_{i}(t) < 1
    &&\forall i \in \ClN, \quad \forall t \in \ClT.
  \end{align}
\end{lem}
\begin{prf}
  See \ref{sec:CD_pi_dot<1_proof} in the appendix. \qed
\end{prf}
\noindent
From \cref{lem:CD_pi_dot<1}, the consistency condition~\eqref{eq:consistency} is satisfied by the candidate queueing delay pattern $\CD{\Vtpi}(t)$.
Based on these observations, we substitute $\CD{\Vtw}(t)$, $\CD{\Vtpi}(t)$, and $\CD{\Vtrho}$ into the \texttt{[DUE-QP]} and exclude the consistency condition~\eqref{eq:consistency} from the problem formulation.
As a result, we obtain the flow determination problem that only involves the flow variables $\Vtq(t)$ and $\Vty(t)$, as $\CD{\Vtw}(t)$, $\CD{\Vtpi}(t)$, and $\CD{\Vtrho}$ are treated as a given parameter.
The problem is formulated as follows:
\begin{align}
  \text{\texttt{[FLOW-LP]}} \quad 
  \min_{[\Vtq, \Vty] \geq \Vt0}. \quad & 
  \FD{Z}(\Vtq, \Vty)  \equiv \UE{Z}(\Vtq, \Vty, \CD{\Vtw}, \CD{\Vtpi}, \CD{\Vtrho}) 
  = 
  Z^{\mathrm{P}}(\Vtq, \Vty)
  + Z^{\mathrm{Q}}(\CD{\Vtw}, \CD{\Vtpi}) 
  -  Z^{\mathrm{D}}(\CD{\Vtw}, \CD{\Vtpi}, \CD{\Vtrho})
  \\
  \text{s.t.} \quad 
  &\int_{t \in \ClT} \Vtq(t) \mathrm{d}t - \VtQ = \Vt0,
  \label{eq:CP_demand_conservation}
  \\
  &\VtA \Vty(t) - \Vtq(t) = \Vt0
  \hspace{7cm}\forall t \in \ClT,
  \label{eq:CP_flow_conservation}
  \\
  & \VtM \overline{\Vttheta}(t)
   - \Vty(t) \geq \Vt0
  \hspace{7cm}\forall t \in \ClT,
  \label{eq:CP_non_queueing}
\end{align}
where $\overline{\Vttheta}(t)$ is defined as follows:
\begin{align}
  &\overline{\Vttheta}(t) \equiv
  \left(
      \CD{\dot{\Vtw}}(t) - \VtA^{\top}_{+} \CD{\dot{\Vtpi}}(t) + \Vt1
  \right)
  &&\forall t \in \ClT.
  \label{eq:dfn_theta}
\end{align}
Here $\VtA^{\top}_{+}$ denotes the matrix obtained from $\VtA^{\top}$ by retaining only the positive entries (i.e., the matrix corresponding to the upstream-end indicator), so that $(\VtA^{\top}_{+}\Vtx)_{(i,j)}=x_{i}$ (the upstream node component) for each link $(i,j)\in\ClL$.
\noindent
The decision variables in this problem are the flow variables $\Vty(t)$ and $\Vtq(t)$.
Consequently, the terms $Z^{\mathrm{Q}}(\CD{\Vtw}, \CD{\Vtpi})$ and $Z^{\mathrm{D}}(\CD{\Vtw}, \CD{\Vtpi}, \CD{\Vtrho})$ in the objective function are constant with respect to $\Vtq(t)$ and $\Vty(t)$.
Unlike the \texttt{[DUE-QP]}, the objective function $\FD{Z}(\Vtq, \Vty)$ in the \texttt{[FLOW-LP]} contains only linear terms and no differential or quadratic terms.
Moreover, $\overline{\Vttheta}(t)$ in the constraints is also constant.
\par
From \cref{lem:LCPtoQP}, if the \texttt{[FLOW-LP]} is feasible and its optimal value is zero, then the candidate cost pattern is an exact equilibrium travel cost pattern, i.e., the GQRP holds.
In addition, in this case, the solution to the \texttt{[FLOW-LP]} provides the equilibrium flow pattern, which satisfies the original DUE conditions.
We have the following theorem and corollary, which are the main results of this section.
\begin{thm}
  \label{thm:QRP}
  The GQRP holds if and only if the \texttt{[FLOW-LP]} is feasible and its optimal value is zero, i.e., $\FD{Z}(\FD{\Vtq}, \FD{\Vty}) = 0$.
\end{thm}
\begin{cor}
  \label{cor:DUE_solution}
  Suppose that the GQRP holds, the following variables are the solution to the \texttt{[DUE-LCP]}:
  \begin{align}
    \left[
     \UE{\Vtq}, 
     \UE{\Vty},
     \UE{\Vtw},
     \UE{\Vtpi},
     \UE{\Vtrho}
   \right]
    = 
    \left[
     \FD{\Vtq}, 
     \FD{\Vty},
     \CD{\Vtw},
     \CD{\Vtpi},
     \CD{\Vtrho}
   \right].
  \end{align}
\end{cor}

\subsubsection{Sufficient condition for the GQRP}
\label{sec:ConstructiveApproach}
The above approach verifies the validity of the GQRP and computes the equilibrium flow pattern by numerically solving \texttt{[FLOW-LP]}.
In the remainder of this section, we present a constructive method that directly recovers the equilibrium flow pattern from the solution to \texttt{[COST-LP-P]}, without solving \texttt{[FLOW-LP]}.
This method applies when the flow pattern $(\CD{\Vtq},\CD{\Vty})$ obtained from \texttt{[COST-LP-P]} is uniquely determined and the schedule delay cost function satisfies a mild condition.
Under these assumptions, the GQRP holds and the solution to \texttt{[DUE-LCP]} can be derived analytically.
\par
First, we substitute $\CD{\Vtw}(t)$, $\CD{\Vtpi}(t)$, and $\CD{\Vtrho}$ obtained by solving the \texttt{[COST-LP-D]} into the \texttt{[DUE-LCP]}, and reduce the conditions~\eqref{eq:departure_time_choice}, \eqref{eq:route_choice}, and \eqref{eq:queueing} to the following simplified conditions:
\begin{align}
  &\Vt0 \leq \overline{\Vtsigma}(t) \perp \Vtq(t) \geq \Vt0
  &&\forall t \in \ClT,
  \label{eq:QRP_cnd_q}
  \\
  &\Vt0 \leq \overline{\Vtlambda}(t) \perp \Vty(t) \geq \Vt0
  &&\forall t \in \ClT,
  \label{eq:QRP_cnd_y}
  \\
  &\Vt0 \leq 
  \left(
    \VtM \overline{\Vttheta}(t)
   - \Vty(t)
  \right) \perp \CD{\Vtw}(t) \geq \Vt0
  &&\forall t \in \ClT,
  \label{eq:QRP_cnd_w}
\end{align}
where $\overline{\Vtsigma}(t)$ and $\overline{\Vtlambda}(t)$ are defined as follows:
\begin{align}
  &\overline{\Vtsigma}(t) \equiv \left(
    \CD{\Vtpi}(t) + \Vts(t) - \CD{\Vtrho}
  \right)
  &&\forall t \in \ClT,
  \label{eq:dfn_sigma}
  \\
  &\overline{\Vtlambda}(t) \equiv \left(
    \CD{\Vtw}(t) - \VtA^{\top} \CD{\Vtpi}(t) + \Vtc
  \right)
  &&\forall t \in \ClT.
  \label{eq:dfn_lambda}
\end{align}
Our goal is to derive $\UE{\Vty}$ and $\UE{\Vtq}$ that satisfy the above conditions, together with the demand conservation condition~\eqref{eq:demand_conservation} and the flow conservation condition~\eqref{eq:flow_conservation}, and to identify when such $\UE{\Vty}$ and $\UE{\Vtq}$ exist.
To achieve this goal with a clear mathematical characterization, we introduce the following additional assumption on the structure of $\CD{\Vty}$, $\CD{\Vtq}$, and $\CD{\Vtw}$ obtained from \texttt{[COST-LP-P]} and \texttt{[COST-LP-D]}.
The assumption is stated as follows:
\begin{asm}
\label{asm:q>0_p>0}
The optimal solutions to \texttt{[COST-LP-P]} and \texttt{[COST-LP-D]}, namely $\CD{\Vtq}$, $\CD{\Vty}$, and $\CD{\Vtw}$, satisfy the following condition:
  \begin{align}
    &\CD{q}_{i}(t) > 0 \quad \text{and} \quad \CD{y}_{ij}(t) > 0
    \quad \Leftrightarrow \quad
    \CD{w}_{ij}(t) > 0
    &&\forall (i,j) \in \ClL,\; t \in \ClT.
  \end{align}
\end{asm}
\cref{asm:q>0_p>0} implies that a candidate queue (i.e., $\CD{w}_{ij}(t) > 0$) can occur on a link only at times when the upstream node generates a positive demand flow.
This assumption corresponds to conditions introduced in prior studies on corridor networks~\citep{Fu2022-st,Sakai2024-wy} that exclude \emph{false bottlenecks}, i.e., bottlenecks at which candidate queues never form.
Moreover, this assumption provides additional information on the complementarity conditions~\eqref{eq:QRP_cnd_q}, \eqref{eq:QRP_cnd_y}, and \eqref{eq:QRP_cnd_w}, and plays a purely mathematical role in enabling the derivation of closed-form expressions for $\UE{\Vty}$ and $\UE{\Vtq}$.\footnote{
\cref{asm:q>0_p>0} also eliminates cases in which arbitrariness exists in the flow patterns $\CD{\Vtq}$ and $\CD{\Vty}$ obtained from \texttt{[COST-LP-P]}.
Here, by \emph{arbitrariness} we mean situations in which only aggregate quantities are determined, while their decomposition remains indeterminate; for example, when only the sum $q_i(t)+q_j(t)$ for adjacent nodes $i$ and $j$ is determined, or when only the total outflow $\sum_{j\in\ClN^{\mathrm{OUT}}(i)} y_{ij}(t)$ is determined but its allocation among outgoing links is not uniquely specified.
Such indeterminacy arises when, at a given time, the cost difference between origins reduces to a constant free-flow travel-time term, making competing origins (or routes) indifferent.
In contrast, once queueing occurs on links that are actually used, this indifference is resolved and the corresponding flow pattern becomes uniquely determined.
Thus, \cref{asm:q>0_p>0} can exclude such indeterminate cases.
}
\par
Under \cref{asm:q>0_p>0}, we now proceed to derive $\UE{\Vty}$ and $\UE{\Vtq}$ that satisfy the conditions~\eqref{eq:QRP_cnd_q}, \eqref{eq:QRP_cnd_y}, and \eqref{eq:QRP_cnd_w}, along with the demand conservation condition~\eqref{eq:demand_conservation} and the flow conservation condition~\eqref{eq:flow_conservation}.
We begin by considering the link flow $\UE{y}_{ij}(t)$ during queueing periods when $\CD{w}_{ij}(t) > 0$.
From \cref{eq:QRP_cnd_w}, we find that when the queueing period for each link $(i,j) \in \ClL$, the link flow $\UE{y}_{ij}(t)$ is given by:
\begin{align}
  &\UE{y}_{ij}(t) = \overline{\theta}_{ij}(t) \mu_{ij} = \left( 
  1 - \CD{\dot{\pi}}_{j}(t)  
  \right) \mu_{ij}
  &&\forall t \in \ClT, \quad \CD{w}_{ij}(t) > 0, \quad \forall (i,j) \in \ClL.
  \label{eq:sufficient_FR_y_queueing}
\end{align}
Since $\overline{\theta}_{ij}(t) = 1 - \CD{\dot{\pi}}_{j}(t)$ incorporates information on changes in travel costs induced by downstream congestion, the link flow $\UE{y}_{ij}(t)$ in \cref{eq:sufficient_FR_y_queueing} reflects the impact of downstream congestion on upstream link flows.
We now turn to the link flow $\UE{y}_{ij}(t)$ during non-queueing periods, i.e., when $\CD{w}_{ij}(t)=0$.
In particular, on the subset of non-queueing times at which additionally $\CD{\dot{\pi}}_{j}(t)=0$, we have $\overline{\theta}_{ij}(t)=1$, which motivates the choice below.
Motivated by this observation, we construct the link flow $\UE{y}_{ij}(t)$ during non-queueing periods as
\begin{align}
  &\UE{y}_{ij}(t) = \CD{y}_{ij}(t)
  &&\forall t \in \ClT, \quad \CD{w}_{ij}(t) = 0, \quad \forall (i,j) \in \ClL,
  \label{eq:sufficient_FR_y_non_queueing}
\end{align}
where $\CD{y}_{ij}(t)$ denotes the link flow obtained from \texttt{[COST-LP-P]}.
Combining \cref{eq:sufficient_FR_y_queueing} and \cref{eq:sufficient_FR_y_non_queueing}, we obtain a unified expression for the equilibrium link flow as follows:
\begin{align}
  &\UE{y}_{ij}(t) = \left( 1 - \CD{\dot{\pi}}_{j}(t) \right)\CD{y}_{ij}(t)
  &&\forall t \in \ClT, \quad \forall (i,j) \in \ClL.
  \label{eq:FR_y_sol}
\end{align}
The link flow $\UE{y}_{ij}(t)$ defined in \cref{eq:FR_y_sol} satisfies the complementarity condition~\cref{eq:QRP_cnd_y}, because $\CD{\dot{\pi}}_{j}(t)<1$ holds from \cref{lem:CD_pi_dot<1} and $\CD{y}_{ij}(t)\le\mu_{ij}$ follows from the non-queueing condition~\eqref{eq:DSO_non_queueing} in \texttt{[COST-LP-P]}.
\par
Based on \cref{eq:FR_y_sol}, the corresponding demand flow $\UE{\Vtq}(t)$ is obtained from the flow conservation condition applied to $\UE{\Vty}(t)$ as
\begin{align}
  \UE{q}_{i}(t)
  &=
      \sum_{j \in \ClN^{\mathrm{OUT}}(i)} \UE{y}_{ij}(t)
      - \sum_{j \in \ClN^{\mathrm{IN}}(i)} \UE{y}_{ji}(t)
  \notag
    \\
    &=
    \sum_{j\in\ClN^{\mathrm{OUT}}(i)}
     \left(1-\CD{\dot{\pi}}_{j}(t)\right)\CD{y}_{ij}(t)
     - \sum_{j\in\ClN^{\mathrm{IN}}(i)}
     \left(1-\CD{\dot{\pi}}_{i}(t)\right)\CD{y}_{ji}(t).
    \label{eq:FR_q_sol}
\end{align}
For the constructed link flow in \cref{eq:FR_y_sol} and the demand flow in \cref{eq:FR_q_sol} to constitute the exact equilibrium flow pattern, these flow patterns must satisfy the following three requirements:
\begin{itemize}
  \item[(i)] demand flow satisfies the demand conservation condition~\eqref{eq:demand_conservation},
  \item[(ii)] demand flow is nonnegative, i.e., $\UE{\Vtq}(t) \geq 0$ for all $t \in \ClT$,
  \item[(iii)] optimal value is zero, i.e., $\FD{Z}(\UE{\Vtq}, \UE{\Vty}) = 0$.
\end{itemize}
\noindent
Under \cref{asm:q>0_p>0}, we can confirm that requirement~(i) always holds, since the cumulative link flow of $\UE{y}_{ij}(t)$ is equal to that of $\CD{y}_{ij}(t)$ for all links $(i,j) \in \ClL$, as shown below:
    \begin{align}
    \int_{t \in \ClT} \UE{y}_{ij}(t) \mathrm{d}t
    &=  \int_{t \in \ClT} \CD{y}_{ij}(t) \left( 1 - \CD{\dot{\pi}}_{j}(t) \right) \mathrm{d} t
    \notag
    \\
    &= \int_{t; \CD{w}_{ij}(t) > 0} \CD{y}_{ij}(t) \left( 1 - \CD{\dot{\pi}}_{j}(t) \right) \mathrm{d} t
       + 
       \int_{t; \CD{w}_{ij}(t) = 0} \CD{y}_{ij}(t) \mathrm{d} t
    \notag
       \\
    &= \int_{t; \CD{w}_{ij}(t) > 0} \CD{y}_{ij}(t) \mathrm{d} t
       - 
       \int_{t; \CD{w}_{ij}(t) > 0}
       \mu_{ij} \CD{\dot{\pi}}_{j}(t) \mathrm{d} t
       + 
       \int_{t; \CD{w}_{ij}(t) = 0} \CD{y}_{ij}(t) \mathrm{d} t
    \notag
       \\
    &= \int_{t \in \ClT} \CD{y}_{ij}(t) \mathrm{d} t.
    \label{eq:cumulative_flow_equivalence}
  \end{align}
The remaining requirements are (ii) and (iii), but in conclusion, these requirements are satisfied when the schedule cost function satisfies certain conditions.
These technical details are presented in the appendix, but only a sketch is provided here.
First, we find that $\CD{q}_{i}(t)=0 \quad \Rightarrow \quad \UE{q}_{i}(t)=0$ holds from \cref{eq:FR_q_sol} under \cref{asm:q>0_p>0} (proof is provided in \ref{sec:QRP_sufficient_proof} in the appendix).
Therefore, we investigate the case of $\CD{q}_{i}(t) > 0$ to confirm the non-negativity of $\UE{q}_{i}(t)$.
In this case, the $\UE{\Vtq}$ in \cref{eq:FR_q_sol} can be expanded as follows.
\begin{align}
  \text{\cref{eq:FR_q_sol}}
  &= \sum_{j \in \ClN^{\mathrm{OUT}}_{y(t)>0}(i)} \mu_{ij} (1 - \CD{\dot{\pi}}_{j}(t))
     - \sum_{j \in \ClN^{\mathrm{IN}}_{y(t)>0}(i)} \CD{y}_{ji}(t) (1 - \CD{\dot{\pi}}_{i}(t))
      && \CD{q}_{i}(t) > 0
   \notag
  \\
  &= \sum_{j\in\ClN^{\mathrm{OUT}}_{y(t)>0}(i)}
  \left(1-\CD{\dot{\pi}}_{j}(t)\right)\mu_{ij}
  -
  \sum_{j\in\ClN^{\mathrm{IN}}_{y(t)>0}(i)}
  \left(1+\dot{s}(t)\right)\CD{y}_{ji}(t)
  && \CD{q}_{i}(t) > 0
  \label{eq:FR_q_expanded}
\end{align}
where $\ClN^{\mathrm{OUT}}_{y(t)>0}(i)$ and $\ClN^{\mathrm{IN}}_{y(t)>0}(i)$ are defined as follows:
  \begin{align}
  &\ClN^{\mathrm{OUT}}_{y(t)>0}(i) \equiv \left\{ j \in \ClN^{\mathrm{OUT}}(i) \mid \CD{y}_{ij}(t) > 0 \right\},
  \quad 
  \ClN^{\mathrm{IN}}_{y(t)>0}(i) \equiv \left\{ j \in \ClN^{\mathrm{IN}}(i) \mid \CD{y}_{ji}(t) > 0 \right\}.
\end{align}
Thus, the non-negativity condition $\UE{q}_{i}(t) \geq 0$ is equivalent to the minimum of the right-hand side of \cref{eq:FR_q_expanded} being nonnegative.
The derivation is provided in Appendix~\ref{sec:QRP_sufficient_proof}, and the resulting sufficient condition for non-negativity is stated in the subsequent theorem in terms of the schedule delay cost function. 
Moreover, under this condition, substituting the derived flow patterns $\FD{\Vty}(t)$ and $\FD{\Vtq}(t)$ into the objective function of the \texttt{[FLOW-LP]} yields an objective value of zero, i.e., requirement (iii) is satisfied. 
These results are summarized in the following theorem.
\begin{thm}
  \label{thm:QRP_sufficient}
  Suppose that \cref{asm:q>0_p>0} holds and that, for every $(i, t) \in \ClN \times \ClT$ such that $\ClN^{\mathrm{IN}}_{y(t)>0}(i) \neq \emptyset$, the schedule delay cost function satisfies
  \begin{align}
  &\dot{s}(t) <
  \dfrac{\sum_{j \in \ClN^{\mathrm{OUT}}_{y(t)>0}(i)} \mu_{ij}}{\sum_{j \in \ClN^{\mathrm{IN}}_{y(t)>0}(i)} \mu_{ji}} - 1.
  \label{eq:sufficient_condition}
  \end{align}
  Then, the GQRP holds and the $\UE{\Vty}$ in \cref{eq:FR_y_sol} and $\UE{\Vtq}$ in \cref{eq:FR_q_sol} are the equilibrium flow solution to the \texttt{[DUE-LCP]}.
\end{thm}
\begin{prf}
  See \ref{sec:QRP_sufficient_proof} in the appendix.
  \qed
\end{prf}
The sufficient condition~\eqref{eq:sufficient_condition} is derived from the non-negativity requirement of the demand flow $\Vtq$ computed from the candidate cost pattern.
More intuitively, it can be interpreted as a restriction on how fast queues are allowed to dissipate, particularly during late-arrival periods.
Specifically, the condition imposes an upper bound on the slope of the schedule delay cost function $s(t)$ during late arrivals.
When the demand flow is positive (i.e., $\CD{q}_{i}(t)>0$), the optimality condition~\eqref{eq:CD_departure} yields $\dot{s}(t)=-\CD{\dot{\pi}}_{i}(t)$, so a steeper schedule delay cost implies a faster decrease in the earliest arrival time.
Moreover, since $\CD{q}_{i}(t)$ depends on the difference between inflows and outflows, $\CD{\dot{\pi}}_{i}(t)$ affects $\CD{q}_{i}(t)$ only when outflows are positive, in which case the queueing-delay slope satisfies $\CD{\dot{w}}_{ij}(t)=\CD{\dot{\pi}}_{i}(t)-\CD{\dot{\pi}}_{j}(t)$.
Thus, \eqref{eq:sufficient_condition} effectively prevents unrealistically rapid queue dissipation by bounding the slope of the candidate earliest arrival time function.
The condition is imposed only at node--time pairs $(i, t)$ with at least one active inflow link; for pairs with $\ClN^{\mathrm{IN}}_{y(t)>0}(i) = \emptyset$, the non-negativity of $\UE{q}_{i}(t)$ holds automatically, so the condition is vacuously satisfied.
{
\par
Note that condition~\eqref{eq:sufficient_condition} is sufficient but not necessary; the GQRP may still hold even when it is violated, as can be verified via the \texttt{[FLOW-LP]} (\cref{thm:QRP}).
When the GQRP does not hold, the reconstructed demand flow becomes negative at certain times, corresponding to backward-bending cumulative curves at bottlenecks.
This occurs because, when the late-arrival penalty is too steep, users prefer to depart earlier and experience queueing rather than incur the high schedule delay cost, so that queues remain longer than the candidate pattern suggests (see also \citet{Fu2022-st} for an analogous discussion in the DUE-DTC setting).
Nevertheless, such violations tend to be localized in time and space, and the candidate pattern still provides a high-quality approximation, as confirmed numerically in \cref{sec:CaseGQRPNotHold}.
Furthermore, the domain of applicability of the GQRP can be expanded by reducing the effective slope of the schedule delay cost through policy instruments such as time-dependent subsidies or charges (\cref{sec:PolicyForSatisfyingGQRP}).
}

\subsection{Summary of the Procedure}
The overall procedure of the queue replacement approach based on \cref{thm:QRP,cor:DUE_solution,thm:QRP_sufficient} is summarized as follows.
\begin{framed}
  \noindent \underline{\textbf{ Queue Replacement Approach }}
  \begin{description}
    \item[Step 1] \textbf{Cost determination}
    \begin{description}
      \item[Step 1-1] Derive the \texttt{[COST-LCP]} and transform it into the \texttt{[COST-LP-P]} and \texttt{[COST-LP-D]}
      \item[Step 1-2] Numerically solve the \texttt{[COST-LP-D]} and obtain $\CD{\Vtw}(t)$, $\CD{\Vtpi}(t)$, and $\CD{\Vtrho}$
      \item[Step 1-3] Numerically solve the \texttt{[COST-LP-P]} to obtain $\CD{\Vtq}(t)$ and $\CD{\Vty}(t)$ 
    \end{description}
    \item[Step 2] \textbf{Flow determination}
    \begin{description} 
      \item[Step 2-1] If \cref{asm:q>0_p>0} and the sufficiency condition \eqref{eq:sufficient_condition} hold, derive the equilibrium flow pattern $\FD{\Vtq}(t)$ and $\FD{\Vty}(t)$ based on \cref{thm:QRP_sufficient} else go to \textbf{Step 2-2}
      \item[Step 2-2] Formulate the \texttt{[FLOW-LCP]} and \texttt{[FLOW-LP]} using the $\CD{\Vtw}(t)$, $\CD{\Vtpi}(t)$, and $\CD{\Vtrho}$
      \item[Step 2-3] Numerically solve the \texttt{[FLOW-LP]} and calculate the optimal value $\FD{Z}(\FD{\Vtq}, \FD{\Vty})$
      \item[Step 2-4] If the \texttt{[FLOW-LP]} is feasible and $\FD{Z}(\FD{\Vtq}, \FD{\Vty}) = 0$, then the GQRP holds, and the solution to the \texttt{[DUE-LCP]} is given by \cref{cor:DUE_solution}.
    \end{description}
  \end{description}
  \vspace{-12pt}
\end{framed}
The cost determination step involves solving the \texttt{[COST-LCP]} to obtain the cost pattern $[\CD{\Vtw}, \CD{\Vtpi}, \CD{\Vtrho}]$.
Additionally, in the flow determination step, we solve the \texttt{[COST-LP-P]} to obtain the flow pattern $[\CD{\Vtq}, \CD{\Vty}]$.
Based on $[\CD{\Vtq}, \CD{\Vty}]$, we can verify whether the GQRP sufficient condition holds.
If the GQRP sufficient condition \eqref{eq:sufficient_condition} holds, we can obtain the equilibrium flow pattern using the constructive approach.
In other words, given the cost pattern $[\CD{\Vtw}, \CD{\Vtpi}, \CD{\Vtrho}]$, it is possible to analytically construct the equilibrium flow pattern.
Thus, the constructive approach is much simpler than the numerical approach, which requires solving the LP again.
Note that, even when using the numerical approach, the computation is significantly easier than directly solving the \texttt{[DUE-LCP]}, as it only requires solving an LP, which is a numerically tractable class of problems.

\section{Discussion}
\label{sec:Discussion}
\subsection{Policy implications}
\label{sec:PolicyImplications}
\subsubsection{First-best policies and the GQRP}
\label{sec:FirstBestPolicies}
The GQRP is valuable not only as a theoretical foundation for efficiently and systematically solving the DUE-RDTC problem, but also as a conceptual bridge that clarifies the relationship between the dynamic system optimal (DSO) state\footnote{
  The DSO state is defined as the traffic state in which the total transportation cost is minimized without any queueing delay.
} achieved by a first best policy and the DUE state.
In certain situations, the cost determination problem \texttt{[COST-LP-P]} can be interpreted as the DSO problem.
This equivalence emerges when the users' passing times at intermediate nodes are uniquely determined by the arrival time at the destination, such as when the DSO flow forms a tree structure at each time in the DSO state.
In such a situation, the DSO problem can be formulated using the Lagrangian-like coordinate system approach, similar to the DUE formulation.\footnote{{An} equivalence between the DSO and the cost determination problems was shown in the case of departure time choice without route choice (DTC), where the network structure is a corridor. See \citet{Fu2022-st,Sakai2024-dt, Sakai2024-wy}.}
Building on this structural insight, we now formally introduce an assumption that enables this equivalence to hold under the RDTC setting and proceed to formulate the DSO problem accordingly.
\par
We first introduce the following assumption:
\begin{asm}
  \label{asm:DSO_flow_node_passing_time}
  In the DSO state, the node-passing time for each node is uniquely determined by the destination arrival time.
\end{asm}
A typical example that satisfies this assumption is a flow pattern in which the DSO flow forms a tree structure at each time in the DSO state because the route from each origin to the destination is unique at each time.
In the following analysis, we adopt this assumption.
Under the assumption, the DSO problem can be formulated in the Lagrangian coordinate system, as follows:
\begin{align}
  \text{\texttt{[DSO-LP]}} \quad
  \min_{[\Vtq, \Vty] \geq \Vt0}
  \quad 
  & \sum_{i \in \ClN} \int_{t \in \ClT} s(t) q_{i}(t) \mathrm{d}t
    +
  \sum_{(i,j) \in \ClL} \int_{t \in \ClT} c_{ij} y_{ij}(t) \mathrm{d}t
  \\
  \text{s.t.} \quad
  & \int_{t \in \ClT} \Vtq(t) \mathrm{d}t = \VtQ,
  &&\quad [\Vtrho]
  \\
  & \VtA \Vty(t) - \Vtq(t) = \Vt0
  &&\forall t \in \ClT \quad [\Vtpi(t)],
  \\
  &\Vty(t) \leq \Vtmu
  &&\forall t \in \ClT \quad [\Vtp(t)],
\end{align}
where the variables in the brackets are the Lagrangian multipliers for the constraints.
The solution to the \texttt{[DSO-LP]} is denoted by $\SO{\Vtx}$.
The optimality conditions of the \texttt{[DSO-LP]} are as follows:
\begin{align}
  &\int_{t \in \ClT} \SO{\Vtq}(t) \mathrm{d} t = \VtQ,
  \label{eq:SO_demand_conservation}
  \\
  &\VtA \SO{\Vty}(t) - \SO{\Vtq}(t) = \Vt0
  &&\forall t \in \ClT,
  \label{eq:SO_flow_conservation}
  \\
  &\Vt0 \leq 
  \left(
  \Vtmu - \SO{\Vty}(t)
  \right) \perp \SO{\Vtw}(t) \geq \Vt0
  &&\forall t \in \ClT,
  \label{eq:SO_queue}
  \\
  &\Vt0 \leq \left(
    \SO{\Vtpi}(t) + \Vts(t) - \SO{\Vtrho}
  \right) \perp \SO{\Vtq}(t) \geq \Vt0
  &&\forall t \in \ClT,
  \label{eq:SO_departure}
  \\
  &\Vt0 \leq 
  \left(
    \SO{\Vtw}(t) - \VtA^{\top} \SO{\Vtpi}(t) + \Vtc
  \right) \perp \SO{\Vty}(t) \geq \Vt0
  &&\forall t \in \ClT.
  \label{eq:SO_route}
\end{align}
It is worth noting that these optimality conditions can be interpreted as the equilibrium conditions under an optimal pricing scheme that eliminates queues at all bottlenecks. In this context, $\SO{p}_{ij}(t)$ represents the optimal price on link $(i,j)$ for a user whose destination arrival time is $t$. Similarly, $\SO{\pi}_{i}(t)$ can be seen as the minimum travel cost at node $i$ for such a user, and $\SO{\rho}_{i}$ as the corresponding equilibrium travel cost. Accordingly, \cref{eq:SO_departure} expresses the equilibrium condition for users' departure time choice under the pricing scheme, while \cref{eq:SO_route} expresses the equilibrium condition for route choice. This equilibrium-based interpretation provides a foundation for deriving the analytical solution to \texttt{[DSO-LP]}.
\footnote{
  Another interpretation of $\SO{p}_{ij}(t)$ is that it represents the market-clearing price pattern under a time-dependent tradable bottleneck permit scheme~\citep{Wada2013-xf,Akamatsu2017-bi}.
  The optimality condition in \cref{eq:SO_queue} can then be interpreted as the demand-supply equilibrium condition in the permit market.
}
\par
Obviously, the optimality condition of the \texttt{[DSO-LP]} has {the} same structure as the optimality condition of the \texttt{[COST-LP-P]} and \texttt{[COST-LP-D]}.
This implies that in settings where the (G)QRP holds, first-best strategies such as dynamic optimal congestion pricing can achieve socially optimal states without reducing the utility of any user.  
From this equivalence and the GQRP, we can derive the following theorem:
\begin{thm}
  \label{thm:DSO_DUE}
  Suppose that the GQRP and \cref{asm:DSO_flow_node_passing_time} hold, then, the DSO cost solution is identical to the DUE cost solution, i.e., the following conditions hold:
  \begin{align}
    &\SO{\Vtp}(t) = \UE{\Vtw}(t)
    &&\forall t \in \ClT,
    \\
    &\SO{\Vtpi}(t) = \UE{\Vtpi}(t)
    &&\forall t \in \ClT,
    \\
    &\SO{\Vtrho} = \UE{\Vtrho}.
  \end{align}
\end{thm}
\begin{prf}
  Under \cref{asm:DSO_flow_node_passing_time}, the DSO problem can be formulated using the Lagrangian-like coordinate system approach, where the state variables are defined by the arrival time at the destination node.
  The DSO problem has {the} same structure as the {\texttt{[COST-LCP]}}.
  Thus, the solution to the DSO problem and {\texttt{[COST-LCP]}} are equivalent.
  Hence, \cref{thm:DSO_DUE} is proved.
  \qed
\end{prf}
Note that the theorem {does} not state that the DSO flow pattern is identical to the DUE flow pattern.
The flow patterns in the DSO and DUE states are generally different, as the DSO state does not have queues, while the DUE state has queues.
\par
We now discuss the queue-eliminating pricing pattern $\SO{\Vtp}(t)$ and its relationship with the equilibrium queueing delay pattern $\UE{\Vtw}(t)$.
From \cref{thm:DSO_DUE}, we find that a Pareto improvement can be achieved if the road manager imposes dynamic pricing that mimics the queueing delay pattern observed in the DUE state.  
This result is supported by two key properties presented in \cref{thm:DSO_DUE}:
\begin{enumerate}
  \item The DSO pricing patterns are identical to the DUE queueing delay patterns at each bottleneck;
  \item The trip cost for each user in the DSO state (schedule delay cost + congestion price) is equal to the total trip cost in the DUE state (schedule delay cost + queueing delay).
\end{enumerate}
\par
Based on these properties, we {conclude} that the road manager's total cost {can} be reduced by implementing a pricing scheme equal to the equilibrium queueing delay.
This is because such a scheme allows the manager to collect revenue equivalent to the total queueing delay cost.
In other words, the manager's benefit is improved without increasing any user's equilibrium trip cost, thereby achieving a Pareto improvement.
This result {is formalized} in the following theorem.
\begin{thm}[Pareto improvement]
  \label{thm:DSO_Pareto_improvement}
  Suppose that the GQRP and \cref{asm:DSO_flow_node_passing_time} hold.
  If the road manager imposes the dynamic pricing equal to the queueing delay at all bottlenecks, the road manager's benefit improves without increasing anyone's equilibrium travel cost.
\end{thm}
\begin{prf}
  Based on \cref{thm:DSO_DUE}, the DSO cost solution is identical to the DUE cost solution under \cref{asm:DSO_flow_node_passing_time}.
  Thus, the commuting cost in the DSO state is equal to the commuting cost in the DUE state.
  In addition, the DSO state does not have queues, thus, the queueing delay cost is zero.
  Therefore, the road manager's benefit is improved by the revenue from the dynamic pricing equal to the queueing delay cost in the DUE state.
  This proves the theorem.
  \qed
\end{prf}
In this way, the GQRP not only offers numerical advantages but also clarifies the relationship between DUE and DSO states in terms of welfare analysis.
If future research succeeds in clarifying the economic interpretation of the cost determination problem \texttt{[COST-LP-P]} in cases where \cref{asm:DSO_flow_node_passing_time} does not hold, it is expected that the GQRP will be shown to provide not only computational benefits but also meaningful economic insights.

\subsubsection{Policy for satisfying the GQRP}
\label{sec:PolicyForSatisfyingGQRP}
The previous sections have shown that the GQRP is effective both as a computational tool and as a foundation for evaluating first-best policies.
This section focuses on the sufficient conditions for the GQRP to hold and discusses policy measures that may help ensure the GQRP.
\par
As shown in \cref{thm:QRP_sufficient}, under \cref{asm:q>0_p>0}, a sufficient condition for the GQRP is given by an explicit inequality involving the bottleneck capacities and the slope of the schedule delay cost function.
This highlights the role of schedule-delay sensitivity in the condition.
Importantly, without altering travelers' intrinsic preferences, traffic demand management policies can modify the effective schedule delay component of the generalized cost by introducing a time-dependent subsidy or charge $\tau(t)$ (e.g., time-of-day incentives).
In such cases, the relevant term becomes $s(t)+\tau(t)$, and its slope can be shaped to facilitate satisfaction of the sufficient condition even when the physical capacity pattern is fixed.
This idea is formalized in the following proposition.
\begin{pro}
  \label{pro:GQRP_scaled_s}
  Consider the following scaled schedule delay cost function:
  \begin{align}
    \tilde{s}(t)=
    \begin{dcases}
      s(t) & \text{if } t<t^{\ast},\\
    \kappa s(t) & \text{if } t\ge t^{\ast},
  \end{dcases}
\end{align}
where $\kappa\in(0,1]$.
Suppose that \cref{asm:q>0_p>0} holds.
Then there exists $\kappa\in(0,1]$ such that the slope restriction in \cref{thm:QRP_sufficient} is satisfied for $\tilde{s}(t)$, and hence the GQRP holds under $\tilde{s}(t)$.
\end{pro}
\begin{prf}
  See \ref{sec:GQRP_scaled_s_proof} in the appendix. \qed
\end{prf}
\noindent
This proposition shows that, for any given capacity pattern, the slope restriction can be satisfied by an appropriate scaling of the schedule delay cost function.
This implies that the validity of the GQRP is not determined by network capacities alone: policies that affect the travelers' effective schedule-delay costs can play an active role in making the GQRP more likely to be satisfied.

\begin{figure}[tbp]
  \centering
    \includegraphics[width=0.85\textwidth]{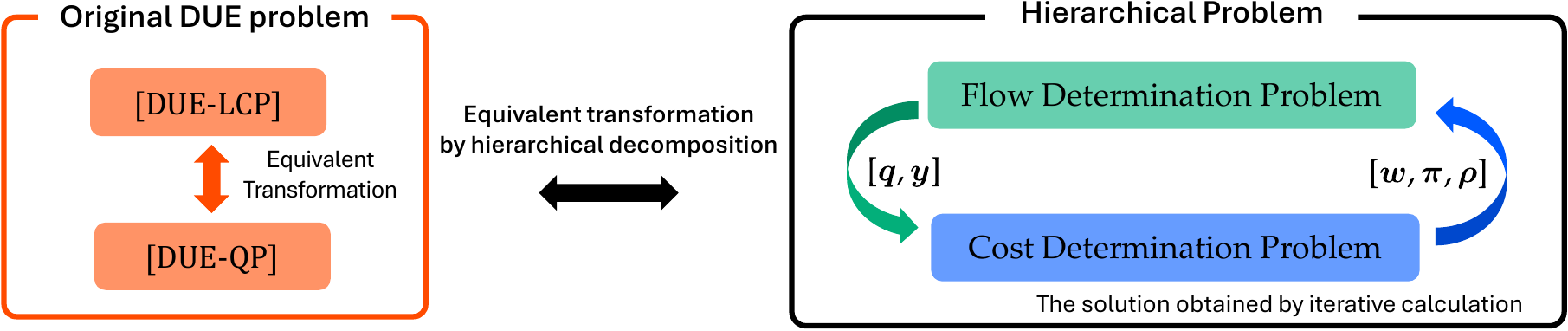}
    \caption{Flowchart of the hierarchical numerical algorithm}
    \label{fig:algo_image}
\end{figure}

\subsection{Algorithmic implications}
\label{sec:AlgorithmicImplications}
The GQRP also provides valuable insights for numerically solving the DUE problem in cases where the principle does not hold. 
The proposed two-step procedure yields exact equilibrium solution when the GQRP holds.
This fact suggests that decomposing the DUE problem into separate subproblems for determining cost and flow variables remains a promising strategy even when the GQRP fails.
This perspective suggests that the original DUE formulation can be viewed as a hierarchical problem: a cost-determination step followed by a flow-determination step, as illustrated in \cref{fig:algo_image}. 
In general, such a hierarchical structure would be solved iteratively by alternating between cost and flow updates until equilibrium is reached.
Mathematically, the cost determination step first produces a candidate cost pattern $[\Vtw, \Vtpi, \Vtrho]$, which is then used in the flow determination step to compute the corresponding flow pattern $[\Vtq, \Vty]$. Based on this updated flow pattern, a new candidate cost pattern $[\Vtw, \Vtpi, \Vtrho]$ is computed, and the process is repeated until convergence. The proposed approach based on the GQRP can thus be interpreted as a special case in which this iterative loop converges in a single step.
\par
Moreover, even when the GQRP does not hold, the cost determination problem \texttt{[COST-LP-D]} is often expected to provide a high-quality initial solution for iterative algorithms.
As discussed earlier, the failure of the GQRP typically occurs when the candidate queueing delay patterns $\CD{\Vtw}$ dissipate too rapidly, violating queueing condition.
Empirically, our numerical results in the latter section indicate that such violations tend to be localized in time or space, i.e., restricted to a few links or specific time periods.
This implies that in most parts of the network and time period, the candidate queueing delay pattern remains consistent with the true equilibrium queueing delay pattern.
Therefore, the solution obtained from the cost determination step is expected to provide an effective initial solution for numerically solving the original DUE formulation, such as \texttt{[DUE-LCP]} or \texttt{[DUE-QP]}.
\par
If the effectiveness of the cost determination solution as a high-quality initial guess can be theoretically established, then the proposed framework, which combines the queue replacement approach with a numerical refinement step, is expected to offer a practical, computationally efficient, and scalable solution method for the DUE-RDTC problem.

\section{Numerical Examples}
\label{sec:IllustrativeExamples}
This section presents numerical experiments to evaluate the effectiveness of the proposed approach.
\cref{sec:BraessNetwork} applies the queue replacement approach to a small-scale Braess network, illustrating step-by-step how it derives an exact DUE solution.
\cref{sec:MandL_Network} evaluates the applicability of the method to medium- and large-scale networks.
\cref{sec:CaseGQRPNotHold} investigates a setting where the GQRP does not hold, showing that the method can still yield a high-quality approximate solution.
\par
The networks used in this section are summarized in \cref{tab:networks}.

\begin{table}[!ht]
  \centering
  \begin{tabular}{|c|c|c|c|}
  \hline
  Network name & \# of nodes & \# of links & \# of origins \\ \hline \hline
  Braess network         & 4   & 5   &  3   \\ \hline
  Nguyen network         & 12   & 18   & 11   \\ \hline
  Sioux Falls network    & 24   & 76   & 23   \\ \hline
  Eastern-Massachusetts network        & 74   & 258   & 73    \\ \hline
  \end{tabular}
  \caption{Summary of networks used in the numerical experiments.}
  \label{tab:networks}
\end{table}

\subsection{Braess Network}
\label{sec:BraessNetwork}
\begin{figure}[tbp]
  \centering
  \begin{minipage}[b]{0.30\textwidth}
    \centering
    \includegraphics[width=0.95\textwidth]{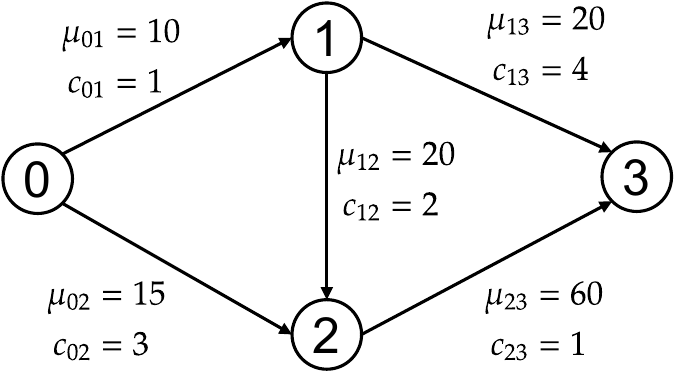}
    \caption{Braess network}
    \label{fig:simple_network}
  \end{minipage}
  \centering
\begin{minipage}[b]{0.34\textwidth}
  \centering
  \includegraphics[width=0.99\textwidth]{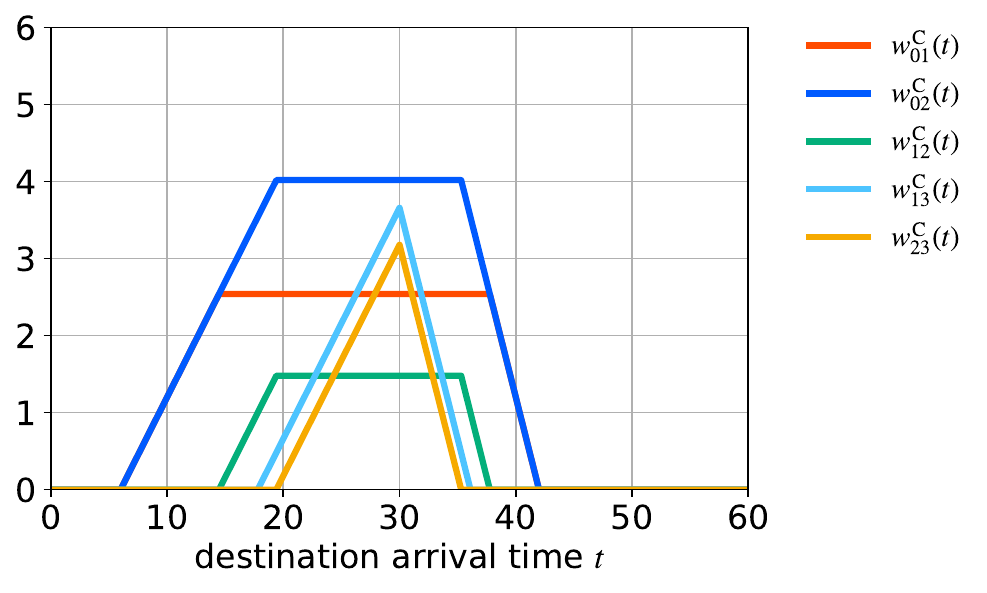}
  \caption{$\CD{\Vtw}(t)$ in Braess network}
  \label{fig:Braess_w}
\end{minipage}
\centering
\begin{minipage}[b]{0.34\textwidth}
  \centering
  \includegraphics[width=0.99\textwidth]{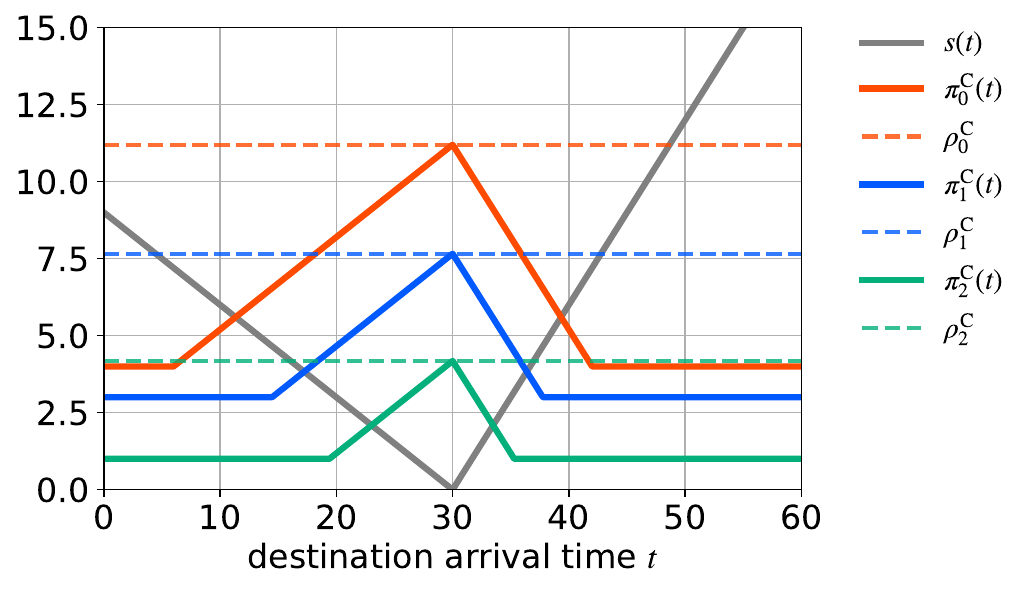}
  \caption{$\CD{\Vtpi}(t)$ and $\CD{\Vtrho}$ in Braess network}
  \label{fig:Braess_pi_rho}
\end{minipage}
\end{figure}

This section provides an illustrative, step-by-step example of the queue replacement approach using a small-scale Braess network consisting of 4 nodes and 5 links, as depicted in \cref{fig:simple_network}.  
The network comprises three origin nodes $\{0, 1, 2\}$ and one destination node $\{3\}$.  
The free flow travel times and bottleneck capacity patterns are shown in \cref{fig:simple_network}.
The OD demand pattern is defined as $Q_{0} = 900$, $Q_{1} = 600$, and $Q_{2} = 400$.  
The schedule delay cost function is specified as follows:
\begin{align}
  &s(t) = 
  \begin{dcases}
    0.3 (t^{\mathrm{P}}-t)
    & \text{if} \quad t \leq t^{\mathrm{P}}
    \\
    0.6(t-t^{\mathrm{P}})
    & \text{if} \quad t > t^{\mathrm{P}}
  \end{dcases}
  \label{eq:numerical_schedule_delay_cost_function}
\end{align}
where the preferred arrival time is $t^{\mathrm{P}} = 30$. 
We numerically solve both problems using the standard simplex method with an arrival time discretization interval of $\Delta t = 0.1$.  
\par
Following the procedure outlined in \textbf{Queue Replacement Approach}, we derive the DUE solution for this simple case.  
In \textbf{Step 1}, we formulate the \texttt{[COST-LCP]} and transform it into its primal and dual LPs, \texttt{[COST-LP-P]} and \texttt{[COST-LP-D]}.  
This yields the solution $\left[\CD{\Vtq}, \CD{\Vty}, \CD{\Vtw}, \CD{\Vtpi}, \CD{\Vtrho} \right]$.
\Cref{fig:Braess_w} shows $\CD{\Vtw}(t)$, which represents the candidate queueing delay pattern.
\Cref{fig:Braess_pi_rho} shows $\CD{\Vtpi}(t)$ and $\CD{\Vtrho}$. 
In both figures, the horizontal axis indicates the arrival time at the destination.
\par
Subsequently, based on \textbf{Step 2-1}, we check the sufficiency condition \eqref{eq:sufficient_condition}.
If the sufficiency condition holds, we can analytically derive the equilibrium flow pattern $\FD{\Vty}(t)$ and $\FD{\Vtq}(t)$ based on \cref{thm:QRP_sufficient}.
If the sufficiency condition does not hold, we proceed to \textbf{Step 2-2} and formulate the \texttt{[FLOW-LP]} using $\CD{\Vtw}(t)$, $\CD{\Vtpi}(t)$, and $\CD{\Vtrho}$.
If the \texttt{[FLOW-LP]} is feasible and its optimal value is zero, then the GQRP holds, and we obtain the equilibrium flow pattern $\FD{\Vty}(t)$ and $\FD{\Vtq}(t)$.
In this settings, the sufficiency condition \eqref{eq:sufficient_condition} holds; thus, we derive the equilibrium flow pattern $\UE{\Vty}(t)$ and $\UE{\Vtq}(t)$ based on \cref{thm:QRP_sufficient}.
Following this procedure, we finally obtain the DUE solution.
\par
\begin{figure}[tbp]
\centering
  \includegraphics[width=0.85\textwidth]{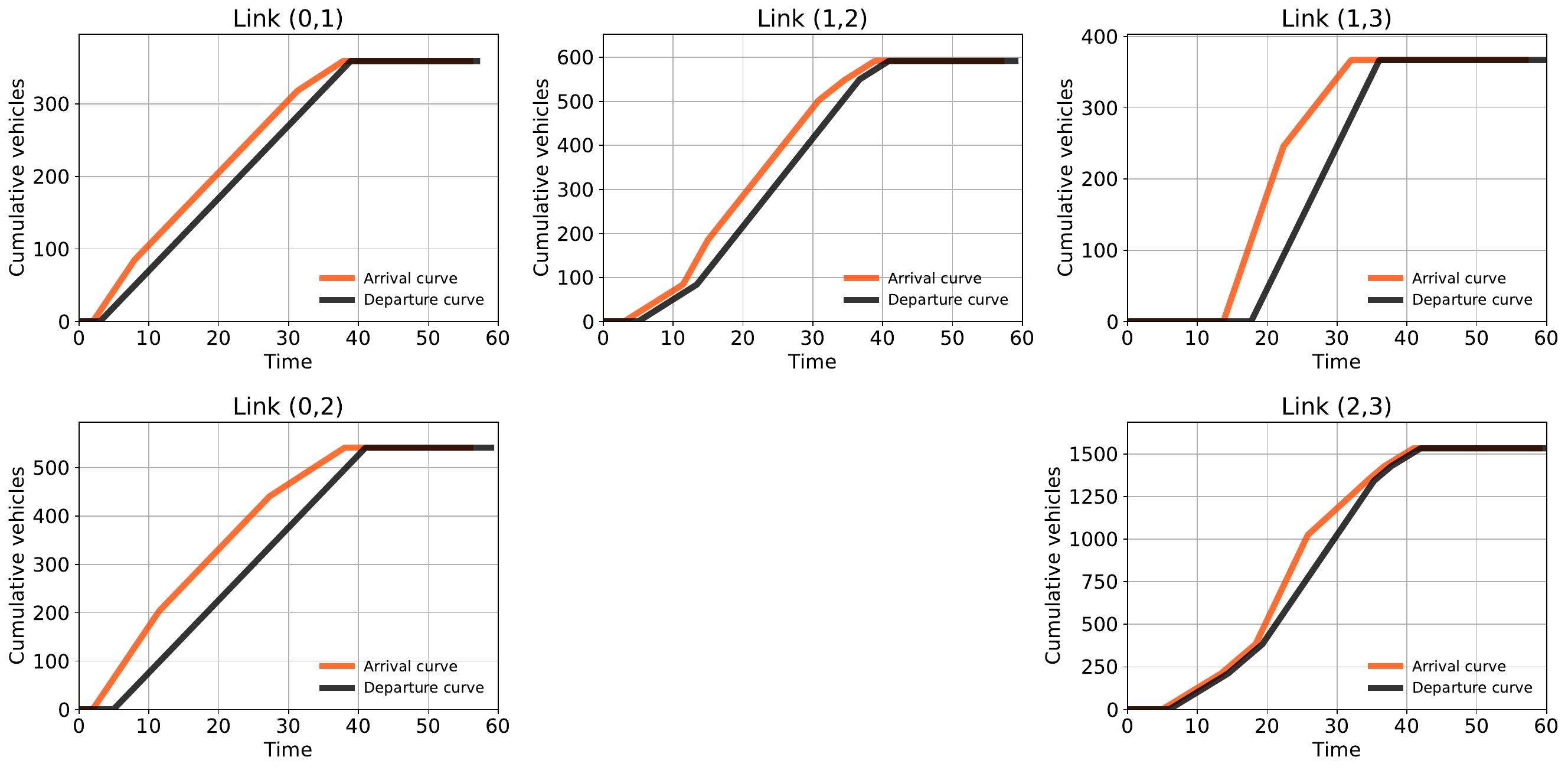}
  \caption{Cumulative equilibrium arrival and departure curves at each link (Braess network)}
  \label{fig:cumulative_flow}
\end{figure}

\cref{fig:cumulative_flow} illustrates the cumulative arrival and departure curves at each link, constructed based on the equilibrium flow pattern $\UE{\Vty}(t)$ and node-passing time pattern $\UE{\Vtpi}(t)$.  
These curves provide a visual representation of the temporal consistency between upstream and downstream flows at each link.


\begin{figure}[tbp]
  \begin{minipage}[b]{0.45\textwidth}
    \centering
    \includegraphics[width=0.6\textwidth]{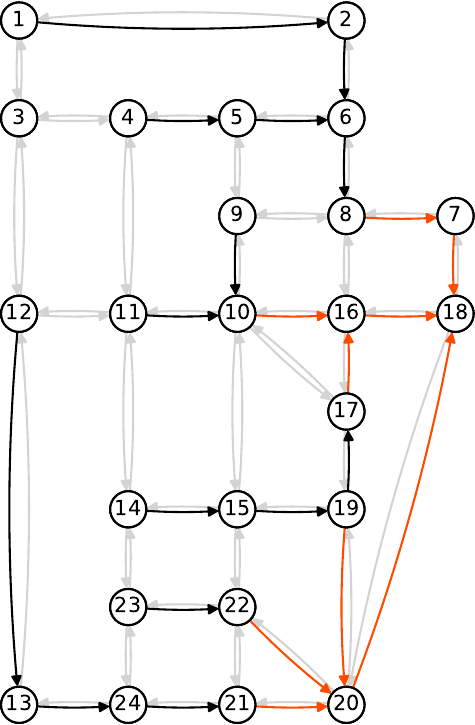}
    \caption{Spatial pattern of equilibrium queueing delay in Sioux Falls network}
    \label{fig:SF_Queue_Network}
  \end{minipage}
  \centering
  \begin{minipage}[b]{0.49\textwidth}
    \centering
    \includegraphics[width=0.85\textwidth]{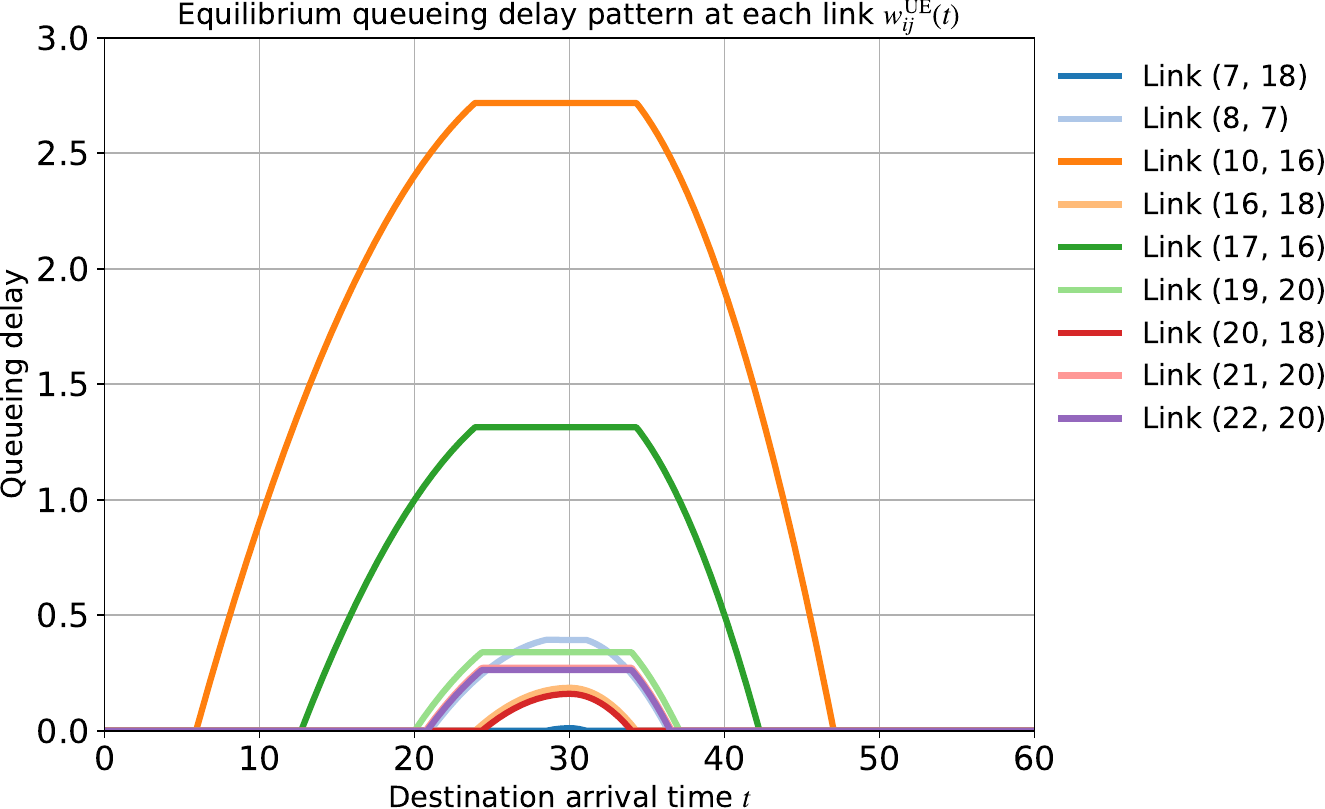}
    \caption{Temporal pattern of equilibrium queueing delay in Sioux Falls network}
    \label{fig:SF_w}
  \end{minipage}
\end{figure}

\begin{figure}[tbp]
  \centering
  \includegraphics[width=0.99\textwidth]{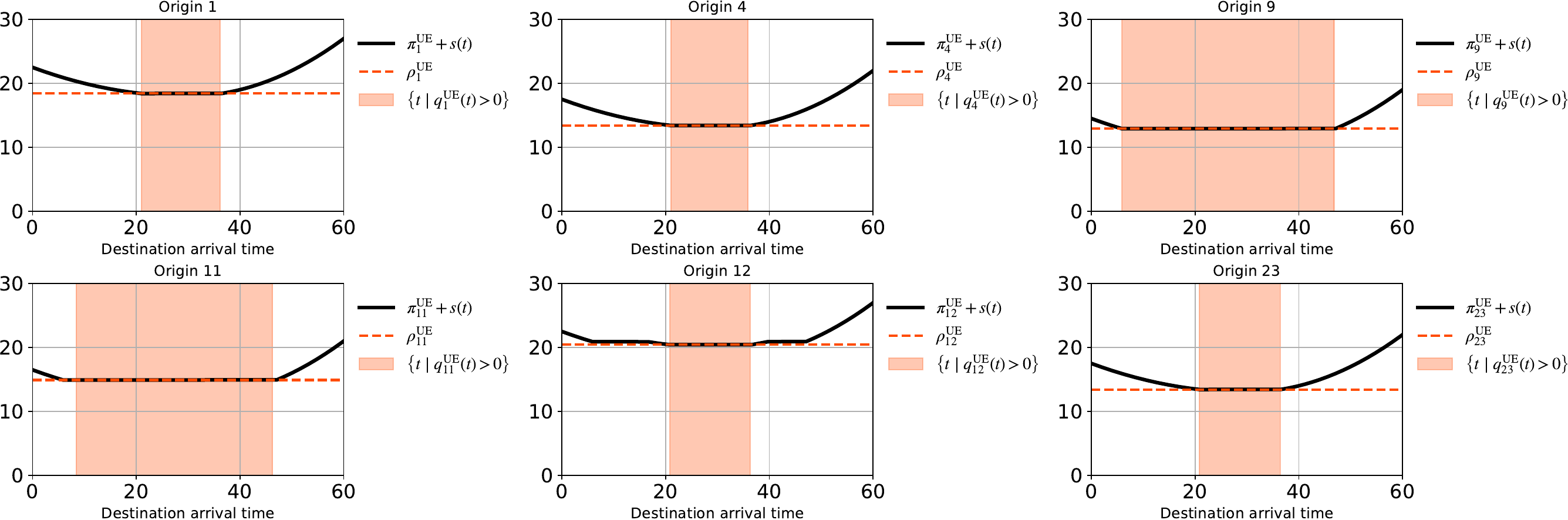}
  \caption{Equilibrium cost and demand flow in the Sioux Falls network for selected origins $1, 4, 9, 11, 12, 23$}
  \label{fig:SF_Cost_DemandFlow}
\end{figure}

\begin{figure}[tbp]
  \centering
    \includegraphics[width=0.95\textwidth]{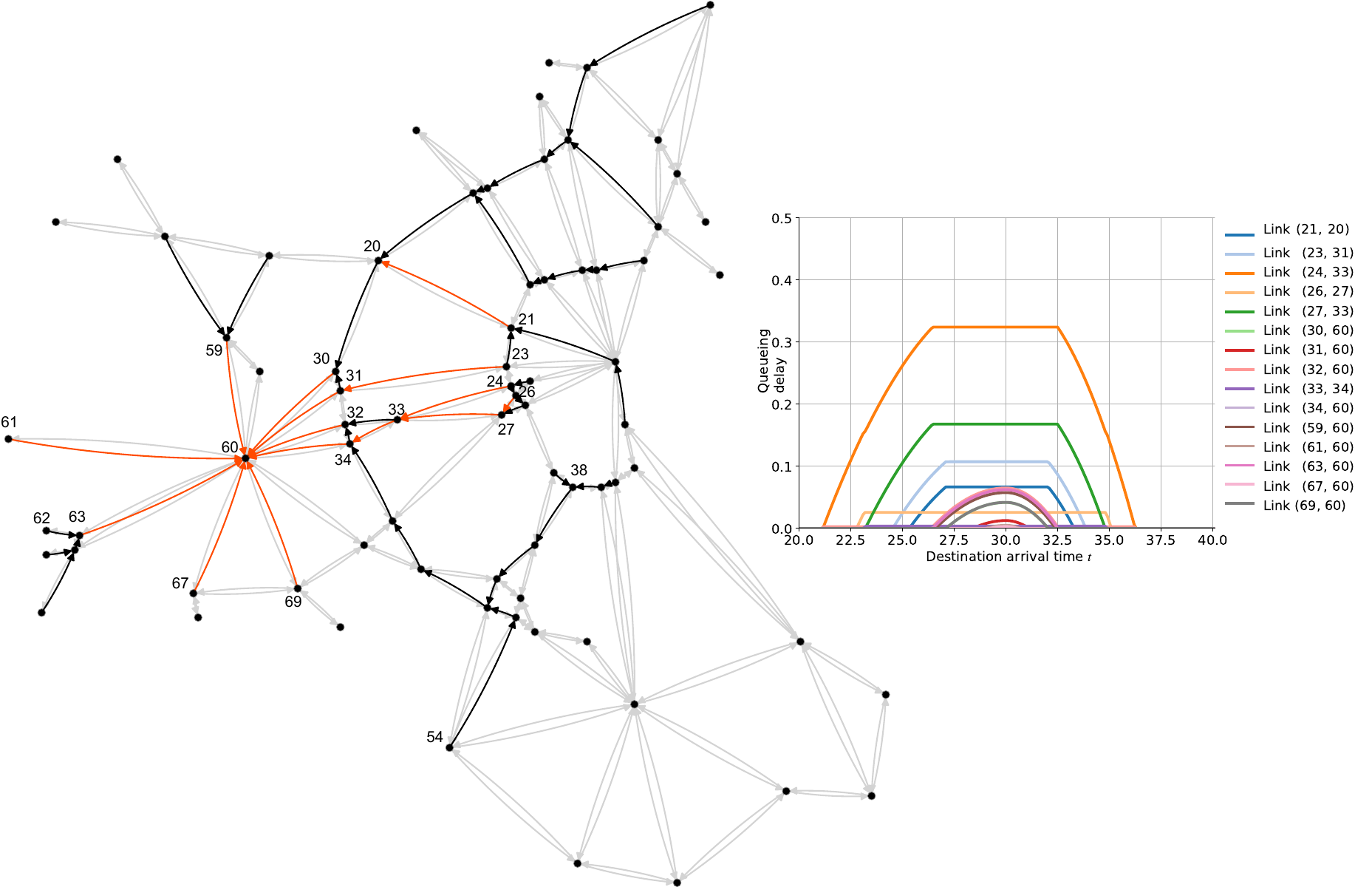}
    \caption{Equilibrium queueing pattern in Eastern-Massachusetts network}
    \label{fig:EMA_Queue_Network}
 \end{figure}

 \begin{figure}[tbp]
    \centering
    \includegraphics[width=0.99\textwidth]{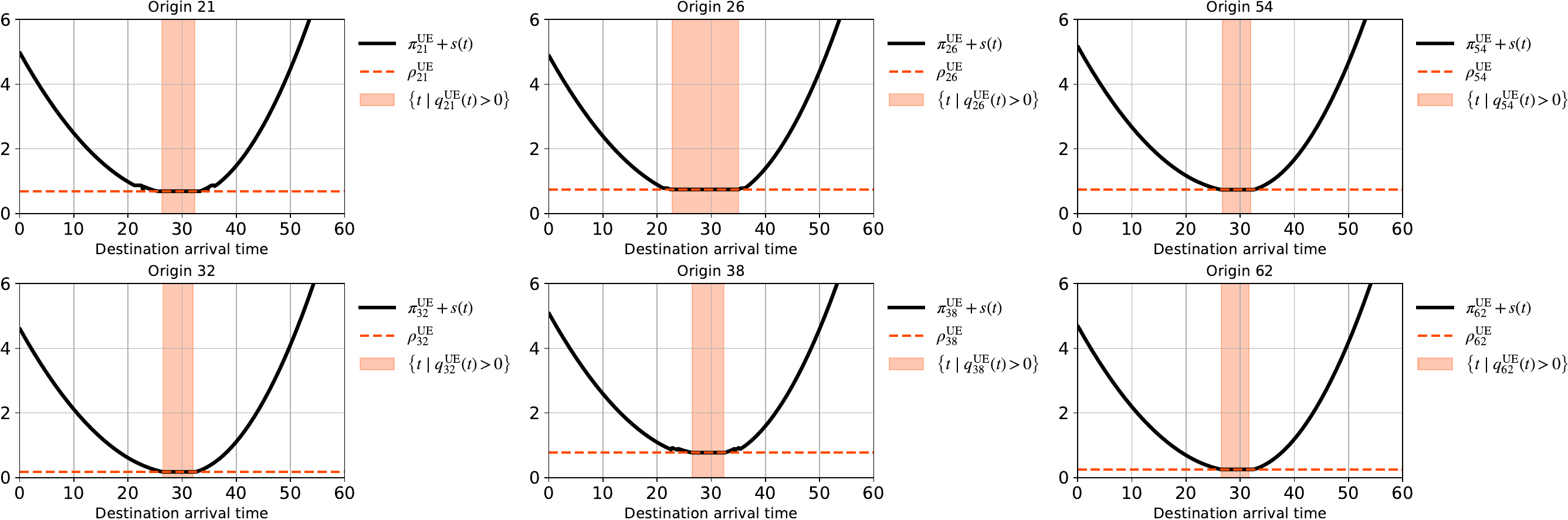}
    \caption{Equilibrium cost and demand flow in Eastern-Massachusetts network for selected origins $21, 26, 54, 32, 38, 62$}
    \label{fig:EMA_Cost_DemandFlow}
 \end{figure}

\subsection{Sioux Falls and Eastern-Massachusetts network}
\label{sec:MandL_Network}
\par
We demonstrate that the proposed queue replacement approach can provide an exact DUE solution even for larger networks.  
We apply the proposed approach to the Sioux Falls network and the Eastern-Massachusetts network, both provided by \citet{transportation_networks_research}~(Accessed: July 22, 2025).  
The numbers of nodes, links, and origins in both networks are summarized in \cref{tab:networks}.  
The dataset includes the capacity pattern, free flow travel times, and the OD demand.
We scale the original capacity pattern by a factor of $1/200$ because the original capacity level is too large for the many-to-one OD pair setting considered in this study.
The schedule delay cost function is specified as follows:
\begin{align}
  s(t) = \begin{dcases}
    \dfrac{\beta}{60} (t^{\mathrm{P}} - t)^{2} & \text{if} \quad t \leq t^{\mathrm{P}} \\
    \dfrac{\gamma}{60} (t - t^{\mathrm{P}})^{2} & \text{if} \quad t > t^{\mathrm{P}}
  \end{dcases}
  \label{eq:parabolic_schedule_delay_cost_function}
\end{align}
where the preferred arrival time is $t^{\mathrm{P}} = 30$, and $\beta, \gamma > 0$ are parameters of the value of the schedule delay.  
We fix $\beta = 0.3$ and $\gamma=0.6$.
Assignment time window is set to $\ClT = [0, 60]$, and the arrival time discretization interval is set to $\Delta t = 0.1$.
\par
We apply the proposed approach to the Sioux Falls and Eastern Massachusetts networks.
For Sioux Falls, node~$18$ is set as the destination.
\cref{fig:SF_Queue_Network} illustrates the spatial distribution of queues: red links indicate queueing links (i.e., links where a queue forms at least once during the assignment period), black links exhibit flow without queues, and gray links have no flow.
\cref{fig:SF_w} shows the temporal pattern of queueing delays for each link, with the horizontal axis representing arrival times at the destination.
\cref{fig:SF_Cost_DemandFlow} presents the equilibrium travel cost and demand flow patterns for a subset of origins ($1, 4, 9, 11, 12, 23$); results for other origins are omitted for brevity.
In each figure, the black line represents the cost pattern $\UE{\pi}_{i}(t) + s(t)$, and the red dashed line denotes the equilibrium travel cost $\UE{\rho}_{i}$.
The shaded red area marks the time window in which queues form, $\left\{ t \in \ClT \mid \UE{q}_{i}(t) > 0\right\}$; within this window, the cost pattern matches $\UE{\rho}_{i}$, confirming that the equilibrium condition~\eqref{eq:departure_time_choice} is satisfied.
Although not shown in detail here, we verified that the equilibrium conditions hold with high accuracy for all origins.
Moreover, the obtained solution achieves an objective value of $\UE{Z}(\cdot) < 10^{-6}$ for the \texttt{[DUE-QP]}, indicating that it is essentially identical to the exact DUE solution.
\par
For the Eastern Massachusetts network, the node~$60$ is set as the destination.
\cref{fig:EMA_Queue_Network} shows the spatial distribution and temporal pattern of queue, using the same color scheme as in Sioux Falls: red for queueing links, black for links with flow but no queues, and gray for links without flow.
\cref{fig:EMA_Cost_DemandFlow} presents the equilibrium travel cost and demand flow patterns for selected origins ($21, 26, 54, 32, 38, 62$); results for other origins are omitted for brevity.
In each case, the black line indicates $\UE{\pi}_{i}(t) + s(t)$, the red dashed line represents $\UE{\rho}_{i}$, and the shaded red area highlights the departure time intervals with positive flow.
Within these intervals, the cost pattern aligns with $\UE{\rho}_{i}$, confirming that the departure time choice condition~\eqref{eq:departure_time_choice} holds.
As with Sioux Falls, we confirmed that the equilibrium conditions are satisfied with high numerical accuracy for all origins.
The solution also achieves an objective value of $\UE{Z}(\cdot) < 10^{-6}$ for the \texttt{[DUE-QP]}, demonstrating that it is effectively equivalent to the exact DUE solution.
\par
In summary, the proposed method successfully produced high-accuracy DUE solutions for both the Sioux Falls and Eastern Massachusetts networks.
In both cases, the equilibrium travel cost patterns matched the theoretical conditions over the active departure time intervals, and the resulting objective function values were less than $10^{-6}$ for \texttt{[DUE-QP]}.
These results demonstrate the scalability and robustness of the approach across networks of different sizes and topology.
{
Regarding the scalability of the proposed approach, the computational cost is dominated by solving two LPs: \texttt{[COST-LP-P]}/\texttt{[COST-LP-D]} and \texttt{[FLOW-LP]}.
The number of decision variables in each LP is $O((N+L)|\ClT|)$, and the number of constraints is also $O((N+L)|\ClT|)$, where $N = |\ClN|$ is the number of origin nodes, $L = |\ClL|$ is the number of links, and $|\ClT|$ is the number of time discretization steps.
Importantly, each LP exhibits the structure of a minimum-cost flow problem with capacity constraints, which is a well-studied problem class for which efficient algorithms and mature solvers are available.
Since the problem size scales linearly with $N$, $L$, and $|\ClT|$, the proposed method is expected to be applicable to networks of considerably larger scale.
}

\subsection{No-GQRP case}
\label{sec:CaseGQRPNotHold}
\begin{figure}[tbp]
  \centering
  \begin{minipage}{0.33\textwidth}
    \centering
    \includegraphics[width=0.80\textwidth]{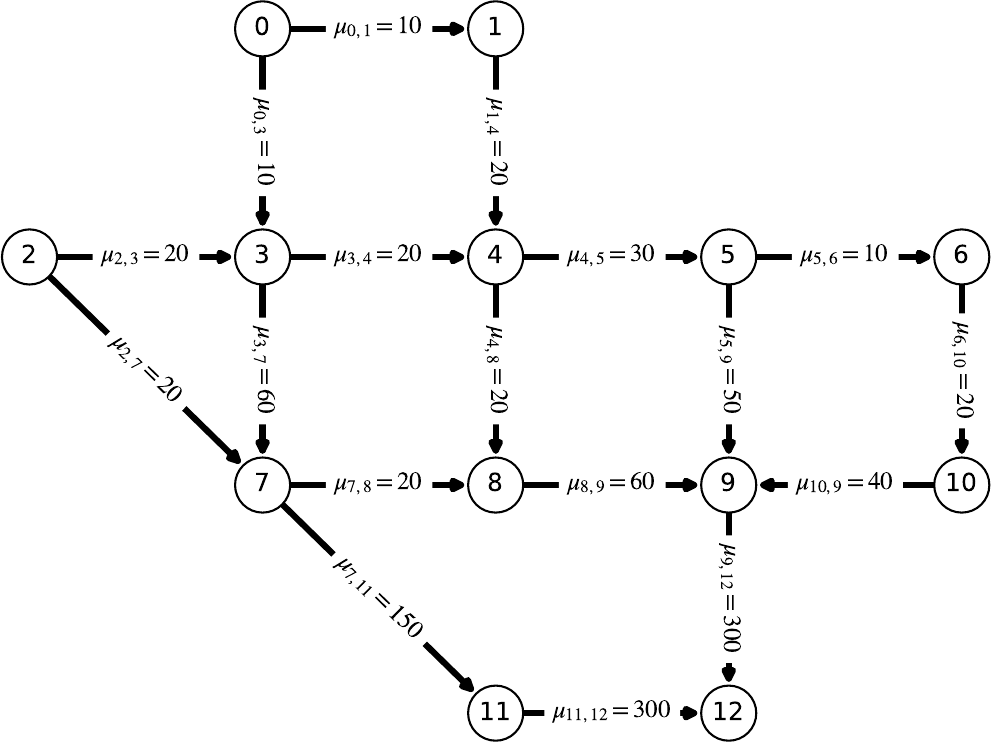}
    \caption{Nguyen network with bottleneck capacity pattern}
    \label{fig:NguyenNetwork}
  \end{minipage}
  \centering
  \begin{minipage}{0.33\textwidth}
    \centering
    \includegraphics[width=0.80\textwidth]{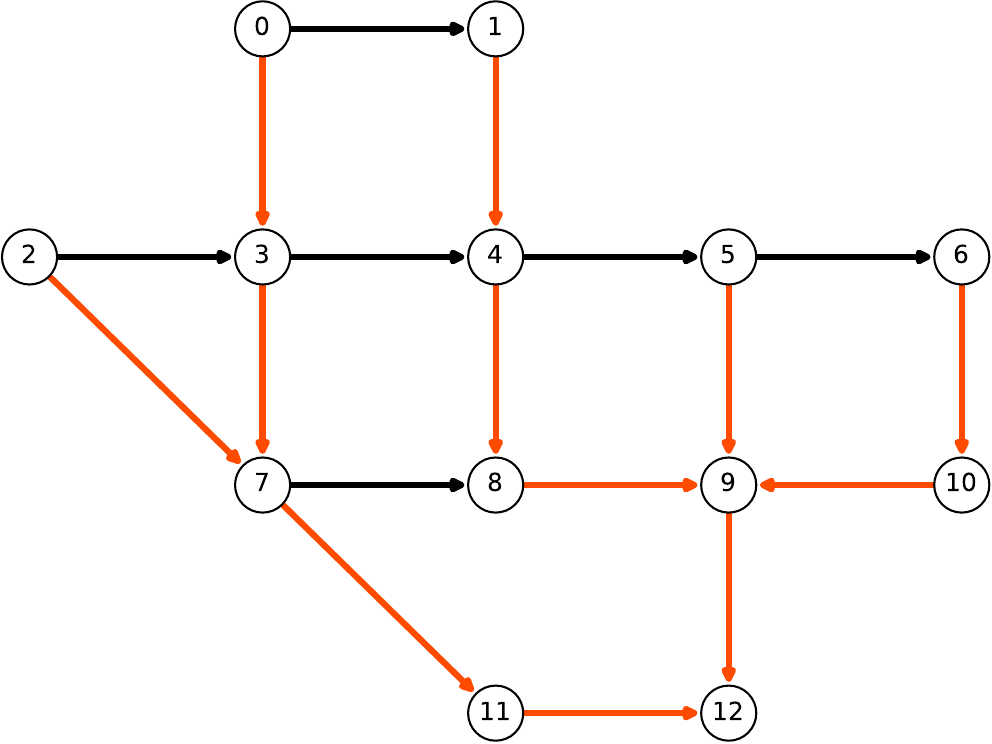}
    \caption{Queueing pattern: GQRP holds}
    \label{fig:Nguyen_GQRP_queue}
  \end{minipage}
  \centering
  \begin{minipage}{0.33\textwidth}
    \centering
    \includegraphics[width=0.80\textwidth]{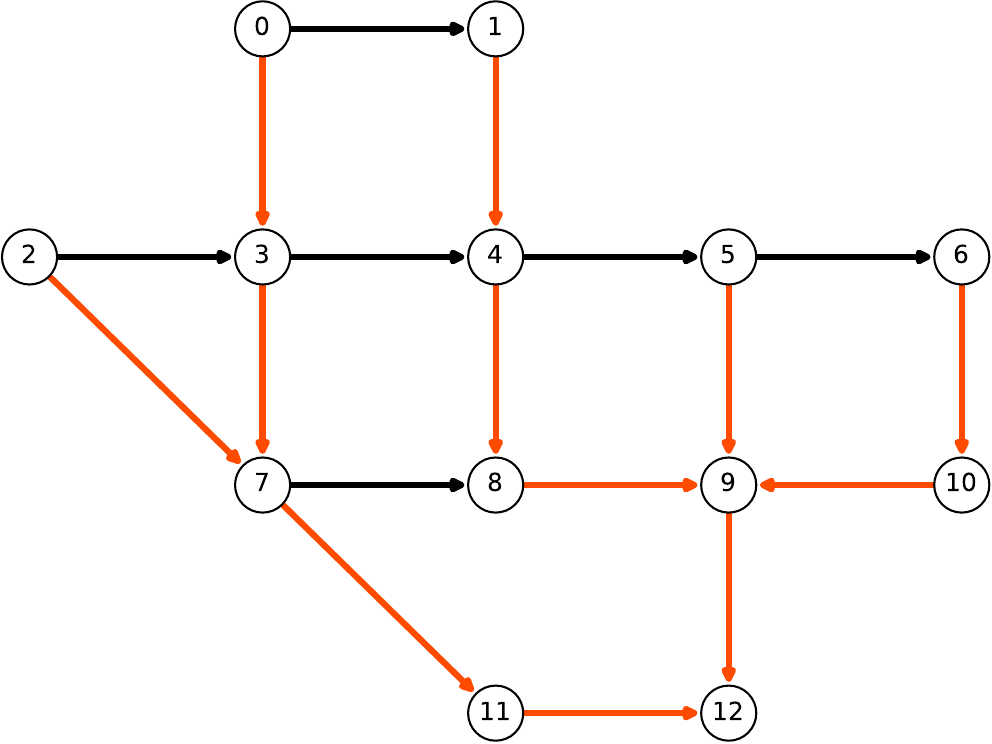}
    \caption{Queueing pattern: GQRP does not hold}
    \label{fig:Nguyen_NoGQRP_queue}
  \end{minipage}
\end{figure}

To further explore the implications of cases where the GQRP does not hold, as discussed in \cref{sec:Discussion}, we conduct a comparative numerical study on the Nguyen network (\cref{fig:NguyenNetwork}), which comprises 12 nodes and 18 links. This experiment is designed to highlight both the validity conditions for the GQRP and the effectiveness from a numerical calculation perspective of the proposed method even when the GQRP does not hold.
\par
We consider two scenarios differentiated by the slope $\gamma$ of the schedule delay cost function $s(t)$ defined in \eqref{eq:numerical_schedule_delay_cost_function}, with preferred arrival time $t^{\mathrm{P}} = 30$ and $\beta = 0.4$. For both scenarios, the capacity pattern is shown in \cref{fig:NguyenNetwork}, and the demand is set to 400 at nodes $0$-$6$ and 300 at nodes $7$-$11$. 
In the GQRP-valid scenario, we set $\gamma = 0.6$; in the no-GQRP scenario, we set $\gamma = 1.2$.
\par
The candidate queueing delay patterns obtained from the cost determination step are shown in \cref{fig:Nguyen_GQRP_queue} and \cref{fig:Nguyen_NoGQRP_queue}.
Interestingly, the spatial distribution of queueing links is identical in both scenarios.
Next, we solve the flow determination problem \texttt{[FLOW-LP]} using these candidate queueing delay patterns.
In the GQRP-valid scenario, the objective function value is $\FD{Z}(\cdot) = 3.00 \times 10^{-11}$, indicating that all DUE conditions are satisfied to high numerical accuracy.
In contrast, in the no-GQRP scenario, the value increases to $\FD{Z}(\cdot) = 86.3$, showing that the obtained solution does not fully satisfy the DUE conditions.
To determine which conditions are violated, we examine the individual components of the objective function.
The analysis reveals that most violations stem from the queueing condition~\eqref{eq:queueing}, with smaller deviations in the departure time choice condition~\eqref{eq:departure_time_choice} and route choice condition~\eqref{eq:route_choice}.
Importantly, these violations are extremely small in magnitude and highly localized—limited to a few links or specific time intervals.
This indicates that even when the GQRP does not hold, the proposed method still produces approximate solutions that are very close to the exact DUE.
Consequently, in GQRP-invalid settings, the solutions from the proposed method are expected to provide high-quality initial points for numerical algorithms solving the original DUE problem (\texttt{[DUE-LCP]} or \texttt{[DUE-QP]}), potentially accelerating convergence substantially.
These findings underscore the practical usefulness and broad applicability of our approach across diverse network structures and user behavior settings.

\section{Conclusion}
\label{sec:Conclusion}
This study presented a hybrid analytical-numerical framework for solving the dynamic user equilibrium (DUE) with route and departure time choice (RDTC) in general many-to-one networks.
The core of the framework is the \textit{Generalized Queue Replacement Principle} (GQRP), which establishes an equivalence between the equilibrium queueing delay pattern and the solution to a linear program (LP) obtained by relaxing the original DUE-RDTC problem.
The GQRP allows the DUE-RDTC problem to be decomposed into two separate LPs.
The first step is the cost determination LP, which provides a candidate cost and queueing delay pattern.
The second step is the flow determination LP, which computes the equilibrium flow pattern based on the candidate queueing delay pattern.
This decomposition makes it possible to obtain the DUE-RDTC solution without iterative or recursive updates of cost and flow variables.
We developed a systematic procedure to verify whether the GQRP holds.
Theoretical analysis identified sufficient conditions for its validity, which depend on the slope of the schedule delay cost function.
We also showed that this slope can be adjusted through policy measures, thereby increasing the likelihood that the GQRP holds.
The GQRP further provides policy and algorithmic insights, offering a clear link between DUE and DSO states.
Numerical experiments on benchmark networks, including Sioux Falls and Eastern Massachusetts, demonstrated that the proposed method achieves high numerical accuracy.
These results indicate that the proposed framework is robust, computationally efficient, and applicable to large-scale traffic networks.
\par
Future research directions include extending the current framework to accommodate heterogeneous user classes, generalizing the OD patterns beyond many-to-one settings, and developing a complete characterization of the necessary and sufficient conditions under which the GQRP holds. 
Moreover, exploring the theoretical relationship between the DSO and DUE states remains an important avenue for further investigation. 
Although this study provides some preliminary findings on this relationship under special conditions in the appendix, a comprehensive and systematic understanding has yet to be established. Clarifying this relationship based on the GQRP would be highly beneficial for welfare analysis and the design of effective traffic management policies. 
These efforts will further strengthen the theoretical foundations and enhance the practical applicability of the DUE-RDTC framework.

\section*{Acknowledgements}

The authors are grateful to the associate editor and anonymous referees for valuable comments on an earlier version of the paper.
We thank Yuki Takayama and Takeshi Nagae for their helpful comments and discussions.
We also thank Masanao Wakui for his assistance in organizing the network data used in the numerical experiments.
This work was supported by Council for Science, Technology and Innovation (CSTI), Cross-ministerial Strategic Innovation Promotion Program (SIP), the 3rd period of SIP ``Smart Infrastructure Management System'' [Grant JPJ012187] (Funding agency: Public Works Research Institute, Japan).
This work was also supported by Japan Society for the Promotion of Science (JSPS) KAKENHI Grant Numbers JP24K22973, JP26K17486, JP24K00999, and JP24K01002, and the Support Program for Urban Studies of the Obayashi Foundation.

\section*{Declaration of generative AI use}

During the preparation of this work the authors used Claude and GPT AI families in order to coding and writing assistance. After using this tool/service, the authors reviewed and edited the content as needed and takes full responsibility for the content of the published article.

\appendix

\section{Proof of Lemma and Theorem}
\subsection{Proof of Lemma \ref{lem:CD_pi_dot<1}}
\label{sec:CD_pi_dot<1_proof}
\begin{figure}[!ht]
  \centering
  \includegraphics[width=0.5\columnwidth]{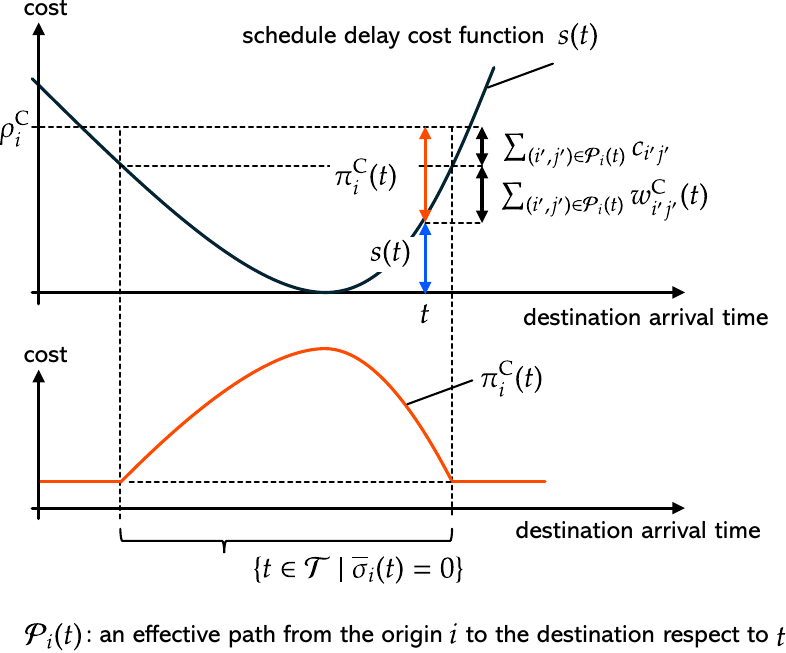}
  \caption{Graphical relationship between $s(t)$ and $\CD{\pi}_{i}(t)$}
  \label{fig:CD_pi_image}
\end{figure}
\begin{prf}
  We consider two cases according to the sign of $\overline{\sigma}_{i}(t) \equiv \CD{\pi}_{i}(t)+s(t)-\CD{\rho}_{i}$.

  \textbf{Case 1}: $\overline{\sigma}_{i}(t)=0$, i.e., $\CD{q}_{i}(t)>0$.
  The departure-choice condition~\eqref{eq:CD_departure} gives $\CD{\pi}_{i}(t)+s(t)=\CD{\rho}_{i}$, and differentiating yields $\CD{\dot{\pi}}_{i}(t)=-\dot{s}(t)$.

  \textbf{Case 2}: $\overline{\sigma}_{i}(t)>0$, i.e., $\CD{q}_{i}(t)=0$.
  By \cref{asm:q>0_p>0}, every active outgoing link $(i,j)$ with $\CD{y}_{ij}(t)>0$ satisfies $\CD{w}_{ij}(t)=0$, hence $\CD{\dot{w}}_{ij}(t)=0$.
  Differentiating the route-choice equality $\CD{w}_{ij}(t)=\CD{\pi}_{i}(t)-\CD{\pi}_{j}(t)-c_{ij}$ then gives $\CD{\dot{\pi}}_{i}(t)=\CD{\dot{\pi}}_{j}(t)$.
  Iterating this relation along the active path from $i$, the recursion terminates either at the destination, where $\CD{\dot{\pi}}_{d}(t)=0$, or at an active origin, where the slope equals $-\dot{s}(t)$ by Case~1.
  Since $s(t)$ is common to all users, the value $-\dot{s}(t)$ is identical at every active origin; the recursion is therefore well-defined and yields $\CD{\dot{\pi}}_{i}(t)\in\{-\dot{s}(t),\,0\}$ (see \cref{fig:CD_pi_image} for an illustration).

  In either case, $\CD{\dot{\pi}}_{i}(t)\in\{-\dot{s}(t),\,0\}$, and since $\dot{s}(t)>-1$, we conclude $\CD{\dot{\pi}}_{i}(t)<1$ for all $i\in\ClN$ and $t\in\ClT$.
  \qed
\end{prf}

\subsection{Proof of \cref{thm:QRP_sufficient}}
\label{sec:QRP_sufficient_proof}
We prove \cref{thm:QRP_sufficient} by showing that, under \cref{asm:q>0_p>0} and the condition \eqref{eq:sufficient_condition}, the flow pattern $\UE{\Vty}(t)$ and $\UE{\Vtq}(t)$ presented in \cref{eq:FR_y_sol,eq:FR_q_sol} are feasible solutions to the \texttt{[FLOW-LP]} and achieve an optimal value of zero.
Specifically, we prove the following two lemmas in order:
\begin{lem}
    \label{lem:FR_y_feasibility}
    Suppose that \cref{asm:q>0_p>0} and the condition \eqref{eq:sufficient_condition} hold.
    The flow pattern $\UE{\Vty}(t)$ and $\UE{\Vtq}(t)$ presented in \cref{eq:FR_y_sol,eq:FR_q_sol} are feasible solutions to the \texttt{[FLOW-LP]}.
  \end{lem}
  \begin{lem}
    \label{lem:FR_y_optimal_value}
    Suppose that \cref{asm:q>0_p>0} and the condition \eqref{eq:sufficient_condition} hold.
    The optimal value of the \texttt{[FLOW-LP]} with the flow pattern $\UE{\Vty}(t)$ and $\UE{\Vtq}(t)$ presented in \cref{eq:FR_y_sol,eq:FR_q_sol} is zero, i.e., $\UE{Z}(\UE{\Vtq}, \UE{\Vty}) = 0$.
  \end{lem}
\subsubsection{Proof of \cref{lem:FR_y_feasibility}}
\begin{prf}
  We verify the five feasibility conditions of the \texttt{[FLOW-LP]} in order.
  \begin{itemize}
    \item[(*)] link flow non-negativity condition
    \item[(**)] demand flow non-negativity condition
    \item[(***)] demand conservation condition~\cref{eq:CP_demand_conservation}
    \item[(****)] flow conservation condition~\cref{eq:CP_flow_conservation}
    \item[(*****)] capacity constraint condition~\cref{eq:CP_non_queueing}
  \end{itemize}
  Under \cref{asm:q>0_p>0}, $\UE{\Vty}(t)$ presented in \cref{eq:FR_y_sol} satisfies the condition (*) and (*****) because $\CD{\Vty}(t) \geq \Vt0$ for all $t \in \ClT$ and $\overline{\Vttheta}(t) \geq \Vt0$ for all $t \in \ClT$ (the latter follows from $\overline{\theta}_{ij}(t)=1-\CD{\dot{\pi}}_{j}(t)>0$ on active links by \cref{lem:CD_pi_dot<1}).
  From the definition of $\UE{\Vtq}(t)$ in \cref{eq:FR_q_sol}, the condition (****) is trivially satisfied.
  Moreover, the condition (***) is also satisfied as discussed in \cref{eq:cumulative_flow_equivalence}.
  \par
  Finally, we verify the condition (**) as follows.
  We first verify that the demand flow $\UE{\Vtq}$ induced by the transformation
\eqref{eq:FR_y_sol} satisfies
\begin{align}
  \CD{q}_{i}(t)=0 \quad \Rightarrow \quad \FD{q}_{i}(t)=0,
  \qquad \forall i\in\ClN,\; t\in\ClT .
  \label{eq:CDq0_implies_FDq0}
\end{align}
Fix $t\in\ClT$ such that $\CD{q}_{i}(t)=0$.
By \cref{asm:q>0_p>0}, for this time $t$ and for every link $(i,j)\in\ClL$ with
$\CD{y}_{ij}(t)>0$, the candidate queueing delay satisfies $\CD{w}_{ij}(t)=0$.
Therefore, \eqref{eq:FR_y_sol} yields
\begin{align}
  \UE{q}_{i}(t)
  &= \sum_{j\in\ClN^{\mathrm{OUT}}(i)} \UE{y}_{ij}(t)
     - \sum_{j\in\ClN^{\mathrm{IN}}(i)} \UE{y}_{ji}(t)
  \notag\\
  &= \sum_{j\in\ClN^{\mathrm{OUT}}(i)} \overline{\theta}_{ij}(t)\,\CD{y}_{ij}(t)
     - \sum_{j\in\ClN^{\mathrm{IN}}(i)} \overline{\theta}_{ji}(t)\,\CD{y}_{ji}(t)
  \notag\\
  &= \sum_{j\in\ClN^{\mathrm{OUT}}(i)}
     \left(1-\CD{\dot{\pi}}_{j}(t)\right)\CD{y}_{ij}(t)
     - \sum_{j\in\ClN^{\mathrm{IN}}(i)}
     \left(1-\CD{\dot{\pi}}_{i}(t)\right)\CD{y}_{ji}(t),
  \label{eq:FDq_expand_noqueue}
\end{align}
where we used $\CD{w}_{ij}(t)=0$ on the support of $\CD{y}_{ij}(t)$ at time $t$.
We now show that $\CD{\dot{\pi}}_{j}(t)=\CD{\dot{\pi}}_{i}(t)$ holds on every active link incident to node $i$.
By \cref{asm:q>0_p>0}, $\CD{w}_{ij}(t)=0$ holds on an interval containing $t$ for every active outgoing link $(i,j)$ (likewise for every active incoming link $(j,i)$).
Differentiating the route-choice complementarity \eqref{eq:CD_route} at equality, namely $\CD{w}_{ij}(t)=\CD{\pi}_{i}(t)-\CD{\pi}_{j}(t)-c_{ij}$, we obtain $\CD{\dot{w}}_{ij}(t)=\CD{\dot{\pi}}_{i}(t)-\CD{\dot{\pi}}_{j}(t)=0$.
Hence there exists a common value $\alpha(t)$ such that $\CD{\dot{\pi}}_{j}(t)=\CD{\dot{\pi}}_{i}(t)\equiv\alpha(t)$ for every $j\in\ClN^{\mathrm{OUT}}(i)$ with $\CD{y}_{ij}(t)>0$ and every $j\in\ClN^{\mathrm{IN}}(i)$ with $\CD{y}_{ji}(t)>0$.
Substituting into \eqref{eq:FDq_expand_noqueue} and using the flow conservation for $\CD{\Vty}$, which under $\CD{q}_{i}(t)=0$ gives $\sum_{j\in\ClN^{\mathrm{OUT}}(i)}\CD{y}_{ij}(t)=\sum_{j\in\ClN^{\mathrm{IN}}(i)}\CD{y}_{ji}(t)$,
\begin{align}
  \UE{q}_{i}(t)
  &= \bigl(1-\alpha(t)\bigr)\!\!\sum_{j\in\ClN^{\mathrm{OUT}}(i)}\!\!\CD{y}_{ij}(t)
   - \bigl(1-\alpha(t)\bigr)\!\!\sum_{j\in\ClN^{\mathrm{IN}}(i)}\!\!\CD{y}_{ji}(t)
   = \bigl(1-\alpha(t)\bigr)\,\CD{q}_{i}(t) = 0,
\end{align}
which proves \eqref{eq:CDq0_implies_FDq0}.
Note that the argument does not require specifying the value of $\alpha(t)$; only the symmetry $\CD{\dot{\pi}}_{j}(t)=\CD{\dot{\pi}}_{i}(t)$ on active links is used.

\par\noindent
It remains to show the nonnegativity of $\UE{q}_{i}(t)$ for times $t$ such that
$\CD{q}_{i}(t)>0$.
For such $t$, define
$\ClN^{\mathrm{OUT}}_{y(t)>0}(i)\equiv\{j\in\ClN^{\mathrm{OUT}}(i)\mid \CD{y}_{ij}(t)>0\}$
and
$\ClN^{\mathrm{IN}}_{y(t)>0}(i)\equiv\{j\in\ClN^{\mathrm{IN}}(i)\mid \CD{y}_{ji}(t)>0\}$.
Then,
\begin{align}
  \UE{q}_{i}(t)
  &= \sum_{j\in\ClN^{\mathrm{OUT}}_{y(t)>0}(i)} \overline{\theta}_{ij}(t)\,\CD{y}_{ij}(t)
     - \sum_{j\in\ClN^{\mathrm{IN}}_{y(t)>0}(i)} \overline{\theta}_{ji}(t)\,\CD{y}_{ji}(t)
  \notag\\
  &= \sum_{j\in\ClN^{\mathrm{OUT}}_{y(t)>0}(i)}
     \left(1-\CD{\dot{\pi}}_{j}(t)\right)\CD{y}_{ij}(t)
     - \sum_{j\in\ClN^{\mathrm{IN}}_{y(t)>0}(i)}
     \left(1-\CD{\dot{\pi}}_{i}(t)\right)\CD{y}_{ji}(t).
  \label{eq:FDq_expand_posq}
\end{align}
Since $\CD{q}_{i}(t)>0$, the optimality condition of \texttt{[COST-LP-P]} and \texttt{[COST-LP-D]} implies $1-\CD{\dot{\pi}}_{i}(t)=1+\dot{s}(t)$.
In addition, for every active outgoing link $(i,j)$ with $\CD{y}_{ij}(t)>0$, \cref{asm:q>0_p>0} combined with $\CD{q}_{i}(t)>0$ gives $\CD{w}_{ij}(t)>0$, and the capacity complementarity in \texttt{[COST-LP-P]} then yields $\CD{y}_{ij}(t)=\mu_{ij}$.
For active incoming links $(j,i)$ with $\CD{y}_{ji}(t)>0$, however, the analogous equality $\CD{y}_{ji}(t)=\mu_{ji}$ does not follow from $\CD{q}_{i}(t)>0$ alone, since it would require $\CD{q}_{j}(t)>0$ at the upstream node $j$, which is not assumed. We therefore retain the inequality $\CD{y}_{ji}(t)\le\mu_{ji}$.
Thus \cref{eq:FDq_expand_posq} becomes
\begin{align}
  \UE{q}_{i}(t)
  &=
  \sum_{j\in\ClN^{\mathrm{OUT}}_{y(t)>0}(i)}
  \left(1-\CD{\dot{\pi}}_{j}(t)\right)\mu_{ij}
  -
  \sum_{j\in\ClN^{\mathrm{IN}}_{y(t)>0}(i)}
  \left(1+\dot{s}(t)\right)\CD{y}_{ji}(t).
  \label{eq:FDq_mu_form}
\end{align}

The right-hand side of \cref{eq:FDq_mu_form} is linear in $(\CD{\dot{\pi}}_{j}(t))_{j\in\ClN^{\mathrm{OUT}}_{y(t)>0}(i)}$, and by \eqref{eq:CD_departure} each component takes one of the two values $\CD{\dot{\pi}}_{j}(t)\in\{-\dot{s}(t),\,0\}$ corresponding to $\CD{q}_{j}(t)>0$ or $\CD{q}_{j}(t)=0$, respectively.
By linearity, the worst case is attained when all active outgoing neighbors share the same value, leading to the following two pure extreme cases:
\begin{itemize}
  \item[(i)] $\CD{\dot{\pi}}_{j}(t)=-\dot{s}(t)$ for all $j\in\ClN^{\mathrm{OUT}}_{y(t)>0}(i)$ (binding when $\dot{s}(t)<0$);
  \item[(ii)] $\CD{\dot{\pi}}_{j}(t)=0$ for all $j\in\ClN^{\mathrm{OUT}}_{y(t)>0}(i)$ (binding when $\dot{s}(t)>0$).
\end{itemize}
If $\UE{q}_{i}(t)\ge 0$ holds in both pure cases, then it holds for every mixed configuration as well.

In case (i), \cref{eq:FDq_mu_form} reduces to
\begin{align}
  \UE{q}_{i}(t)
  &=
  \left(1+\dot{s}(t)\right)
  \left(
    \sum_{j\in\ClN^{\mathrm{OUT}}_{y(t)>0}(i)}\mu_{ij}
    -
    \sum_{j\in\ClN^{\mathrm{IN}}_{y(t)>0}(i)}\CD{y}_{ji}(t)
  \right).
\end{align}
By the flow conservation \eqref{eq:CD_flow_conservation} for $\CD{\Vty}$ together with $\CD{y}_{ij}(t)=\mu_{ij}$ on active outgoing links,
\[
  \sum_{j\in\ClN^{\mathrm{OUT}}_{y(t)>0}(i)}\mu_{ij}
  - \sum_{j\in\ClN^{\mathrm{IN}}_{y(t)>0}(i)}\CD{y}_{ji}(t)
  = \CD{q}_{i}(t) > 0.
\]
Hence $\UE{q}_{i}(t)=(1+\dot{s}(t))\CD{q}_{i}(t) > 0$ since $\dot{s}(t)>-1$. Note that this case does not require \eqref{eq:sufficient_condition}.

In case (ii), \cref{eq:FDq_mu_form} yields
\begin{align}
  \UE{q}_{i}(t)
  &=
  \sum_{j\in\ClN^{\mathrm{OUT}}_{y(t)>0}(i)}\mu_{ij}
  -
  \left(1+\dot{s}(t)\right)\sum_{j\in\ClN^{\mathrm{IN}}_{y(t)>0}(i)}\CD{y}_{ji}(t).
  \label{eq:FDq_case_ii}
\end{align}
We split the analysis into two subcases according to whether $\ClN^{\mathrm{IN}}_{y(t)>0}(i)$ is empty.

\textbf{Subcase (ii-a)}: $\ClN^{\mathrm{IN}}_{y(t)>0}(i)=\emptyset$.
The second sum in \eqref{eq:FDq_case_ii} is taken over the empty set and therefore equals $0$.
Moreover, since $\CD{q}_{i}(t)>0$, the flow conservation condition~\eqref{eq:CD_flow_conservation} for $\CD{\Vty}$ implies $\sum_{j\in\ClN^{\mathrm{OUT}}_{y(t)>0}(i)}\CD{y}_{ij}(t) > 0$, hence $\ClN^{\mathrm{OUT}}_{y(t)>0}(i)\neq\emptyset$.
Consequently,
\[
  \UE{q}_{i}(t) = \sum_{j\in\ClN^{\mathrm{OUT}}_{y(t)>0}(i)}\mu_{ij} > 0,
\]
without invoking the condition~\eqref{eq:sufficient_condition}.

\textbf{Subcase (ii-b)}: $\ClN^{\mathrm{IN}}_{y(t)>0}(i)\neq\emptyset$.
In this subcase $\sum_{j\in\ClN^{\mathrm{IN}}_{y(t)>0}(i)}\mu_{ji}>0$, so the right-hand side of \eqref{eq:sufficient_condition} is well-defined and all subsequent manipulations are free of zero-division.
Since $\CD{y}_{ji}(t)\le\mu_{ji}$ and $1+\dot{s}(t)>0$, we bound
\[
  \UE{q}_{i}(t)\ge
  \sum_{j\in\ClN^{\mathrm{OUT}}_{y(t)>0}(i)}\mu_{ij}
  - \left(1+\dot{s}(t)\right)\sum_{j\in\ClN^{\mathrm{IN}}_{y(t)>0}(i)}\mu_{ji}.
\]
Multiplying the condition \eqref{eq:sufficient_condition} for this $(i,t)$ by the positive denominator $\sum_{j\in\ClN^{\mathrm{IN}}_{y(t)>0}(i)}\mu_{ji}$ yields
\begin{align}
  \left(1+\dot{s}(t)\right)\sum_{j\in\ClN^{\mathrm{IN}}_{y(t)>0}(i)}\mu_{ji}
  <
  \sum_{j\in\ClN^{\mathrm{OUT}}_{y(t)>0}(i)}\mu_{ij},
\end{align}
so the lower bound above is strictly positive, and therefore $\UE{q}_{i}(t)>0$.
Combining subcases (ii-a) and (ii-b), $\UE{q}_{i}(t)\ge 0$ holds for all $i\in\ClN$ and $t\in\ClT$.
Summarizing the above arguments, we conclude that the condition (**) is satisfied.
\par
Thus, the flow pattern $\UE{\Vty}(t)$ and $\UE{\Vtq}(t)$ presented in \cref{eq:FR_y_sol,eq:FR_q_sol} are feasible solutions to the \texttt{[FLOW-LP]} and \cref{lem:FR_y_feasibility} is proved.
\qed
\end{prf}

\subsubsection{Proof of \cref{lem:FR_y_optimal_value}}
\begin{prf}
  By substituting the flow patterns $\UE{\Vty}(t)$ and $\UE{\Vtq}(t)$ presented in \cref{eq:FR_y_sol,eq:FR_q_sol} into the objective function of the \texttt{[FLOW-LP]}, we have:
  \begin{align}
    &\FD{Z}(\UE{\Vtq}, \UE{\Vty}) 
    \\
    &= 
    \int_{t \in \ClT} \overline{\Vtlambda}(t) \UE{\Vty}(t) \mathrm{d}t
    + \int_{t \in \ClT} \overline{\Vtsigma}(t) \UE{\Vtq}(t) \mathrm{d}t
    + \int_{t \in \ClT} \left(
      \VtM \overline{\Vttheta}(t) - \UE{\Vty}(t)
  \right) \CD{\Vtw}(t) \mathrm{d}t
  \\
  &=
  \int_{t \in \ClT} \overline{\Vtlambda}(t) \UE{\Vty}(t) \mathrm{d}t
  + \int_{t \in \ClT} \overline{\Vtsigma}(t) \UE{\Vtq}(t) \mathrm{d}t
  + \int_{t \in \ClT} \left(
    \VtM \overline{\Vttheta}(t) - \overline{\Vttheta}(t) \times \CD{\Vty}(t)
  \right) \CD{\Vtw}(t) \mathrm{d}t
  \\
  &=
  \int_{t \in \ClT} \overline{\Vtlambda}(t) \UE{\Vty}(t) \mathrm{d}t
  + \int_{t \in \ClT} \overline{\Vtsigma}(t) \UE{\Vtq}(t) \mathrm{d}t
  + \int_{t \in \ClT} \left( \overline{\Vttheta}(t) \left( 
    \Vtmu - \CD{\Vty}(t)
  \right) \right) \CD{\Vtw}(t) \mathrm{d}t
  \\
  &=
  \int_{t \in \ClT} \overline{\Vtlambda}(t) \UE{\Vty}(t) \mathrm{d}t
  + \int_{t \in \ClT} \overline{\Vtsigma}(t) \UE{\Vtq}(t) \mathrm{d}t
  \\
  &=\int_{t \in \ClT}  \left(
    \CD{\Vtw}(t) - \VtA^{\top} \CD{\Vtpi}(t) + \Vtc
  \right) \overline{\Vttheta}(t) \CD{\Vty}(t) \mathrm{d}t
  +
  \int_{t \in \ClT}
  \left(
      \CD{\Vtpi}(t) + \Vts(t) - \CD{\Vtrho}
    \right)\UE{\Vtq}(t) \mathrm{d}t
  \\
  &=0.
  \end{align}
  The first integrand vanishes because the route-choice complementarity \eqref{eq:CD_route} gives $\bigl(\CD{w}_{ij}(t)-\CD{\pi}_{i}(t)+\CD{\pi}_{j}(t)+c_{ij}\bigr)\CD{y}_{ij}(t)=0$ on every link.
  The second integrand vanishes because the departure-choice complementarity \eqref{eq:CD_departure} gives $\CD{\pi}_{i}(t)+s(t)-\CD{\rho}_{i}=0$ whenever $\CD{q}_{i}(t)>0$, while $\UE{q}_{i}(t)=0$ whenever $\CD{q}_{i}(t)=0$ (established in the proof of \cref{lem:FR_y_feasibility}).
  Thus, \cref{lem:FR_y_optimal_value} is proved.
  \qed
\end{prf}

\subsection{Proof of \cref{pro:GQRP_scaled_s}}
\label{sec:GQRP_scaled_s_proof}
\begin{prf}
Let $\mathcal{S}$ denote the class of schedule delay cost functions satisfying the conditions described in \cref{sec:NetworkAndUser}.
Fix an arbitrary $s(\cdot)\in\mathcal{S}$ and consider the scaled function $\tilde{s}(t)$ defined in \cref{pro:GQRP_scaled_s}.
By construction, $\tilde{s}(\cdot)\in\mathcal{S}$ for any $\kappa\in(0,1]$, and for $t>t^{\ast}$ we have $\dot{\tilde{s}}(t)=\kappa \dot{s}(t)$.
Let $\dot{s}_{\max} \equiv \max \dot{s}(t)$, which is finite.
For each $(i,t)\in\ClN\times\ClT$ with $\ClN^{\mathrm{IN}}_{y(t)>0}(i)\neq\emptyset$, define
\[
  B_{i,t}\equiv\frac{\sum_{j\in\ClN^{\mathrm{OUT}}_{y(t)>0}(i)}\mu_{ij}}{\sum_{j\in\ClN^{\mathrm{IN}}_{y(t)>0}(i)}\mu_{ji}}-1,
\]
which is well-defined since the denominator is positive; the multiplicative form of \cref{eq:sufficient_condition} on this $(i,t)$ is then equivalent to $\dot{s}(t)<B_{i,t}$.
Let $\bar{B}\equiv\inf_{(i,t)} B_{i,t}$ taken over such $(i,t)$, and assume $\bar{B}>0$ (otherwise no $\kappa\in(0,1]$ can ensure the condition).
Choose $\kappa < \min\left\{1,\; \frac{\bar{B}}{\dot{s}_{\max}}\right\}$.
Then $\dot{\tilde{s}}(t)=\kappa \dot{s}(t)< \bar{B}\le B_{i,t}$ holds for all $t\ge t^{\ast}$ and all $(i,t)$ with $\ClN^{\mathrm{IN}}_{y(t)>0}(i)\neq\emptyset$, i.e., $\tilde{s}(\cdot)$ satisfies \cref{eq:sufficient_condition} on these $(i,t)$; for $(i,t)$ with $\ClN^{\mathrm{IN}}_{y(t)>0}(i)=\emptyset$, the condition is vacuous.
Therefore, by \cref{thm:QRP_sufficient} (under \cref{asm:q>0_p>0}), the GQRP holds under the scaled schedule delay cost function $\tilde{s}(\cdot)$.
\qed
\end{prf}

\end{document}